\documentclass[a4paper, 10pt]{article}
\usepackage{subfiles}
\usepackage[toc,page]{appendix}
\usepackage[bottom]{footmisc}
\usepackage{microtype}
\usepackage{tocbibind}
\usepackage{stmaryrd}   
\usepackage{amsmath}
\usepackage{verbatim}
\usepackage{mathrsfs}
\usepackage{mathtools}
\usepackage{amsthm} 
\usepackage{amssymb} 
\usepackage{faktor}
\usepackage[colorlinks={true},linkcolor={blue},citecolor={blue}, urlcolor={blue} ]{hyperref}
\usepackage{braket}
\usepackage{transparent}
\usepackage{graphicx}
\usepackage{color}
\usepackage{rotating}
\usepackage{tikz}
\usepackage{tikz-cd}
\usepackage{authblk}
\usepackage[margin=2.5 cm]{geometry}

\usetikzlibrary{intersections}
\usetikzlibrary{spy} 
\usetikzlibrary{decorations.pathmorphing}
\usetikzlibrary{patterns}
\usetikzlibrary{backgrounds}

\DeclareMathOperator{\Div}{Div}
\DeclareMathOperator{\im}{im}
\DeclareMathOperator{\cha}{char}
\DeclareMathOperator{\spec}{Spec}

\DeclareMathOperator{\hei}{ht}

\DeclareMathOperator{\fr}{Frac}

\DeclareMathOperator{\res}{res}

\DeclareMathOperator{\Pic}{Pic}
\DeclareMathOperator{\len}{length}

\DeclareMathOperator{\Princ}{Princ}

\DeclareMathOperator{\CH}{CH}

\DeclareMathOperator{\Coh}{Coh}

\numberwithin{equation}{section}

\theoremstyle{plain}
\newtheorem{theorem}{Theorem}[section]
\newtheorem{lemma}[theorem]{Lemma}
\newtheorem{proposition}[theorem]{Proposition}

\newtheorem{corollary}[theorem]{Corollary}
\theoremstyle{definition}
\newtheorem{definition}[theorem]{Definition}

\theoremstyle{remark}
\newtheorem{example}[theorem]{Example}
\newtheorem{remark}[theorem]{Remark}

\newcommand{\sbt}{\,\begin{picture}(-1,1)(-1,-3)\circle*{3}\end{picture}\ }

\newcommand{\catname}[1]{\mathbf{#1}}

\date{}

\makeatletter
\newcommand{\address}[1]{\gdef\@address{#1}}
\newcommand{\email}[1]{\gdef\@email{\url{#1}}}
\newcommand{\@endstuff}{\par\vspace{\baselineskip}\noindent\small
\begin{tabular}{@{}l}\scshape\@address\\\textit{E-mail address:} \@email\end{tabular}}
\AtEndDocument{\@endstuff}
\makeatother
  
\begin{document}
 \title{Adelic geometry on arithmetic surfaces I: idelic and adelic interpretation of the Deligne pairing}
\author{Paolo Dolce}
\address{University of Udine, Italy}
\email{paolo.dolce@uniud.it}
\maketitle
\begin{abstract}
For an arithmetic surface $X\to B=\spec O_K$ the Deligne pairing $\left <\,,\,\right >\colon \Pic(X)\times\Pic(X)\to \Pic(B)$ gives the ``schematic contribution'' to the Arakelov intersection number. We present an idelic and  adelic interpretation of the Deligne pairing; this is the first crucial step for a full idelic and adelic interpretation of the Arakelov intersection number. 

For the idelic approach we show that the Deligne pairing can be lifted to a pairing $\left<\,,\,\right>_i:\ker(d^1_\times)\times \ker(d^1_\times)\to\Pic(B)
$, where $\ker(d^1_\times)$ is an important subspace of the two dimensional idelic group $\mathbf A_X^\times$.  On the other hand, the argument for the adelic interpretation is entirely cohomological.
\end{abstract}

\makeatletter
\@starttoc{toc}
\makeatother

\setcounter{section}{-1}

\section{Introduction}

\subsection{Background}
Problems of topological nature on $(\mathbb Q, |\cdot|)$ (where $|\cdot|$ is the usual euclidean absolute value), are commonly solved after an embedding of $\mathbb Q$ in  its completion $\mathbb R$ with respect to $|\cdot|$. In this way we can take advantage of the completeness properties of $\mathbb R$ and the density of $\mathbb Q$ inside $\mathbb R$. Any other  absolute value which doesn't give the standard euclidean topology on $\mathbb Q$ is equivalent (in the sense of absolute values, see footnote) to the $p$-adic absolute value  $|\cdot|_p$ for any prime $p$, so it makes sense to embed $\mathbb Q$ densely in its completion $\mathbb Q_p$ with respect to $|\cdot|_p$. The above way of reasoning can be easily generalized for any number field $K$, and adeles were introduced in 1930s by Chevalley in order to consider simultaneously all the completions of $K$ with respect all possible places\footnote{A place $\mathfrak p$ is an equivalence class of absolute values on $K$ where two absolute values are declared equivalent if they generate the same topology.}. It is not very useful to study simply the product $\prod_\mathfrak p K_\mathfrak p$ of all completions because the resulting space is ``too big'' and fails to be a locally compact additive group. For any non-archimedean place $\mathfrak p$ let $\mathcal O_{\mathfrak p}$ be the closed unit ball in $K_{\mathfrak p}$ under the standard representative of the place $\mathfrak p$, then the ring of adeles is defined as a subset of $\prod_\mathfrak p K_\mathfrak p$, namely:
$$\mathbf A_K:=\sideset{}{'}\prod_\mathfrak p K_\mathfrak p$$
where $\sideset{}{'}\prod$ is the restricted products with respect to the additive  subgroups $\{\mathcal O_\mathfrak p\colon \mathfrak p\, \text{is non-archimedean}\}$. Note that for archimedean places, the unit ball is not an additive subgroup. The most important features of $\mathbf A_K$ were well described in  \cite{tatethesis} and consist mainly in the fact that $\mathbf A_K$ is a locally compact additive group (so it admits a Haar measure), $K$ is discrete in $\mathbf A_K$, and the quotient $\mathbf A_K/K$ is compact. Moreover the Pontryagin dual $\widehat{\mathbf A_K}$ has a very simple description and $\mathbf A_K\cong\widehat{\mathbf A_K}$. 

The multiplicative version of the adelic theory is the idelic theory, and the group of ideles attached to $K$ is defined as:
$$\mathbf A^\times_K=\sideset{}{'}\prod_\mathfrak p K^\times_\mathfrak p$$
where the restricted product is taken with respect to the subgroups $\mathcal O^\times_\mathfrak p:=\{x\in\mathcal O_{\mathfrak p}\colon |x|_\mathfrak p=1\}$. 
In this case $\mathcal O^\times_\mathfrak p$ has a group structure also in the archimedean case. But a number field $K$ can be seen as the function field of the nonsingular arithmetic curve $B=\spec O_K$, where $O_K$ is the ring of integers of $K$,  and we know  that there is a bijection between points of the completed curve $\widehat B$ in the sense of Arakelov geometry and places of $K$. Therefore the ring of adeles attached to $K$ can be described in a more geometric way  related to $\widehat B$:
$$\mathbf A_{\widehat B}:= \sideset{}{'}\prod_{b\in \widehat B}K_b=\mathbf A_K$$
where $K_b$ is still the local field attached to the ``point'' $b$. So, classical adelic theory can be deduced from $1$-dimensional arithmetic geometry. We can adopt a similar approach but starting from $1$-dimensional algebraic geometry: fix a nonsingular algebraic projective curve $X$ over a perfect field $k$ with function field denoted by $k(X)$; then to each closed point $x\in X$ we can associate a non-archimedean local field $K_x$ with its valuation ring denoted by $\mathcal O_x$. The ring of adeles associated to $X$ is then:
$$
\mathbf A_{X}:= \sideset{}{'}\prod_{x\in X}K_x\,;
$$
In this case $\mathbf A_{X}$ is not a locally compact additive group unless $k$ is a finite field. Each $K_x$ is endowed with a structure of locally linearly compact $k$-vector space (in the sense of \cite{lef}), therefore $\mathbf A_X$ is again locally linearly compact and one can show similarly to the arithmetic case that: $k(X)$ is discrete in $\mathbf A_{X}$, the quotient $\mathbf A_X/k(X)$ is a linearly compact $k$-vector space and $\mathbf A_X$ is a self dual $k$-vector space. In other words, from a topological point of view, the passage from arithmetic theory to algebraic theory implies that we substitute the theory of compactness of groups with the theory of linear compactness of vector spaces. In both arithmetic and geometric $1$-dimensional case, adelic and idelic theory give a generalization of the intersection theory (i.e. the theory of degree of divisors):
\begin{itemize}
\item[\sbt] Ideles can be easily seen as a generalization of line bundles  (resp. Arakelov line bundles), so it is natural to give an extension of the theory of divisors (resp. Arakelov divisors) from an idelic point of view.
\item[\sbt] For an algebraic curve $X$ and any divisor $D\in \Div(X)$ we can define a subspace $\mathbf A_X(D)\subset \mathbf A_X$ and a complex $\mathcal A_X(D)$. The cohomology of $\mathcal A_X(D)$ is equal to the usual Zariski cohomology $H^i(X,\mathscr O_X(D))$, therefore we can give an interpretation of $\deg(D)$ in terms of the characteristic of $\mathcal A_X(D)$ which will be called the adelic characteristic. For a completed  arithmetic curve $\widehat{B}$ we cannot define a complex $\mathcal A_{\widehat B}(\widehat D)$ associated to an Arakelov divisor $\widehat D$, since for archimedean points closed unit balls are not additive groups. However, one can recover the Arakelov degree of $\widehat D$ as the product of volumes (with respect to  an opportune choice of Haar measures) of certain closed balls in $K_b$. 
\end{itemize} 
The above theory remains valid also in the case of singular curves, because by normalization we can always reduce to the nonsingular case.

The first attempt to construct a $2$-dimensional adelic/idelic theory from $2$-dimensional geometry was partially made by Parshin in \cite{pa1}, but he treated only rational adeles (i.e. a subset of the actual ring of adeles) for algebraic surfaces. A more structured approach to the adelic theory for algebraic surfaces was given in \cite{pa2}, still with a few mistakes\footnote{for example in \cite{pa2} the definition of the object $K_x$, and consequently of the subspace $A_{02}$, is wrong (compare with section \ref{sect1.2} for more details).}, and it can be summarized in the following way: let's fix a nonsingular, projective surface $(X,\mathscr O_X)$ over a perfect field $k$, then to each ``flag''  $x\in y$ made of a closed point $x$ inside an integral curve $y\subset X$ we can associate the ring $K_{x,y}$ which will be a $2$-dimensional local field if $y$ is nonsingular at $x$, or a finite product of $2$-dimensional local fields if we have a singularity. Roughly speaking a $2$-dimensional local field is a local field whose residue field is again a local field (see subsection \ref{subs1.1}), and in our case $K_{x,y}$ carries two distinct levels of discrete valuations: there is the discrete valuation associated to the containment $x\in y$ and the discrete valuation associated to $y\subset X$. Formally, $K_{x,y}$ is obtained through a process of successive localisations and completions starting from $\mathscr O_{X,x}$. With the symbol $\mathcal O_{x,y}$ we denote the product of valuation rings inside $K_{x,y}$. Similarly to the $1$-dimensional theory, we perform a ``double restricted product'', first over all points ranging on a fixed curve,  and then over all curves in $X$, in order to obtain the ring of adeles for surfaces:
$$\mathbf A_X:=\sideset{}{''}\prod_{\substack {x\in y\\ y\subset X}} K_{x,y}\subset\prod_{\substack {x\in y\\ y\subset X}} K_{x,y} \,.$$
The topology on $K_{x,y}$ can be defined  thanks to the construction by completions and localisations, and by starting with the standard $\mathfrak m_x$-adic topology on $\mathscr O_{X,x}$, then the topology on $\mathbf A_X$ can be defined canonically. The idelic group attached to $X$ is $\mathbf A^\times_{X}$.

 For $2$-dimensional local fields like $K_{x,y}$ there is a well known theory of differential forms and residues (see for example \cite{yek1}); one can globalise the constructions in order to obtain a $k$-character  $\xi^{\omega}:\mathbf A_X\to k$ associated to a rational differential form $\omega\in\Omega^2_{k(X)|k}$ and the differential pairing:
\begin{eqnarray*}
d_\omega:\mathbf A_X\times\mathbf A_X &\to & k\\
(\alpha,\beta) &\mapsto& \xi^\omega(\alpha\beta)\,.
\end{eqnarray*} 
In \cite{fe0} it is shown that: $\xi^{\omega}$ induces the self duality in the category of $k$-vector spaces,  of  $\mathbf A_X$, the subspace $\mathbf A_X/k(X)^{\perp}$ is linearly compact (orthogonal spaces are calculated with respect to $d_\omega$) and $k(X)$ is discrete in $\mathbf A_X$. 

Both $\mathbf A_X$ and $\mathbf A^\times_{X}$ carry some important subspaces which in turn lead to the construction of certain complexes $\mathcal A_X(D)$, for a divisor $D$, and $\mathcal A_X^\times$ (respectively the ``adelic complex'' and the ``idelic complex''). The cohomology of such complexes can be calculated by geometric methods thanks to the following important results:
\begin{equation}\label{Iso_i_coho}
H^i(\mathcal A^\times_X)\cong H^i(X,\mathscr O^\times_X)\,,
\end{equation}
\begin{equation}\label{Iso_a_coho}
H^i(\mathcal A_X(D))\cong H^i(X,\mathscr O_X(D))\,.
\end{equation}
For a proof of (\ref{Iso_i_coho}) and (\ref{Iso_a_coho}) see respectively \cite{dolce_phd} and \cite{fe0}.  Again, idelic and adelic theory give an extension of the intersection theory on $X$:
\begin{itemize}
\item[\sbt] 
In \cite{pa2} it is shown that the group $\Div(X)$ can be lifted to a subspace if $\mathbf A^\times_{X}$ and the intersection pairing on $\Div(X)$ can be extended at the level of ideles.

\item[\sbt] In  \cite{fe0} it is shown that the characteristic of the complex $\mathcal A_X(D)$ can be used to redefine the intersection pairing between two divisors in terms of adeles, even without using isomorphism (\ref{Iso_a_coho}). Such a theory gives an alternative approach to the Riemann-Roch theorem for algebraic surfaces. 
\end{itemize}
\noindent
Most of the $2$-dimensional constructions outlined above are true for any $2$-dimensional, Noetherian, regular scheme, so in particular for arithmetic surfaces. However,  the adelic theory associated to an arithmetic surface $X\to B=\spec O_K$ ($K$ is a number field) is more complicated and less developed. Locally, the rings $K_{x,y}$ have a completely  different structure  between each other as $2$-dimensional local fields, depending whether $y$ is horizontal or vertical. Moreover, there was the global issue of interpreting the archimedean data of the completed surface $\widehat X$, in the sense of Arakelov geometry, in an adelic and idelic way. A first definition of the ``full'' (or completed) ring of adeles $\mathbf A_{\widehat X}$ has been given only recently in \cite{fe1}.

The purpose of this series of papers is to study the adelic and idelic  theory for a completed arithmetic surface and eventually obtain a $2$-dimensional generalisation of the Tate thesis.

\begin{remark}
In \cite{bei} Beilinson shortly described how to attach a $n$-dimensional adelic theory to any $n$-dimensional Noetherian scheme  in a very abstract functorial way. Reworks and clarifications of this approach are \cite{huber} and \cite[8]{MM1}. In particular in \cite[8.4 and 8.5]{MM1} it is proved  that for dimensions $1$ and $2$ our explicit theory agrees with Beilinson theory of adeles.
\end{remark}

\subsection{Results in this paper}
In this first paper we explain the ``schematic part'' of the idelic and adelic lift of the Arakelov intersection number  (i.e. ignoring the fibres at infinity), whereas a full account of the theory will be published subsequently together with other results. 

Given two Arakelov divisors $\widehat D=D+\sum_\sigma \alpha_\sigma X_\sigma $ and $\widehat E=E+\sum_\sigma \beta_\sigma X_\sigma $, where $D,E\in\Div(X)$, one piece of the Arakelov intersection number $\widehat D.\widehat E$ is obtained thanks to the Deligne pairing $\left<\mathscr O_X(D),\mathscr O_X(E)\right>\in \Pic(B)$. We define $d^1_{\times}$ to be an arrow of the idelic complex $\mathcal A^\times_ X$ associated to $X$, then  we construct an idelic Deligne pairing:
$$
\left<\,,\,\right>_i:\ker(d^1_\times)\times \ker(d^1_\times)\to\Pic(B)
$$
which descends to the Deligne pairing through the composition:
$$\ker(d^1_\times)\times \ker(d^1_\times)\rightarrow\Div(X)\times \Div(X)\to \Pic(X)\times\Pic(X)\,.$$
This will be the arithmetic version  of Parshin idelic lift  given in \cite{pa2}. Our approach will be from local to global:  the main idea consists in globalising Kato's local symbol for $2$-dimensional local fields containing a local field (see  \cite{kato} or \cite{otherLiu}), which is the generalisation of the usual tame symbol for valuation fields (see appendix \ref{k_th}). In the arithmetic framework given by the arithmetic surface $\varphi:X\to B$, we have the following constructions: for any point $x$ sitting on a curve $y\subset X$ we define a ring, which is a finite sum of $2$-dimensional local fields, denoted by $K_{x,y}$; moreover $K_b$ is the local field associated to the point $b\in B$ such that $\varphi(x)=b$. Then Kato's symbol translates into a skew symmetric, bilinear map:
$$(\,,\,)_{x,y}:K^\times_{x,y}\times K^\times_{x,y}\to K^\times_b\,.$$
Roughly speaking, by composing it with the valuation $v_b$ on $K_b$ and by summing over all  flags $x\in y$ such that $\varphi(x)=b$, we show that we obtain a well defined integer $n_b$. By repeating the argument for each $b\in B$ we obtain a divisor $\sum_{b\in B} n_b[b]$. At this point we prove that such a pairing descends to the Deligne pairing.

The adelic theory is very similar to the geometric case and the crucial point consists in  considering the arithmetic analogue of the Euler-Poincare characteristic of coherent sheaves, i.e. the determinant of cohomology. We use the cohomological properties of the adelic complex of the base scheme $B$ in order to give the definition of the adelic determinant of cohomology. Then it is enough to use the fact that the Deligne pairing can be expressed in terms of  the (adelic) determinant of cohomology.

\paragraph{Overview of the contents.} Section \ref{sect1} is a quick review of adelic geometry for arithmetic surfaces, where just by simplicity, we ignore the contribution of fibres at infinity. A more comprehensive introduction to adelic geometry can be found in \cite{dolce_phd}. In section \ref{sect2} and \ref{sect3} we construct respectively the idelic and adelic Deligne pairing. Finally, appendix \ref{k_th} is just a collection of the basic notions of algebraic $K$-theory needed in this paper and appendix \ref{det_coho} is a review of the main features of the determinant of cohomology.  

\paragraph{Basic notations.}
All rings are considered commutative and unitary. When we pick a point $x$ in a scheme $X$ we generally mean a \emph{closed point} if not otherwise specified, also all sums $\sum_{x\in X}$ are meant to be ``over all closed points of $X$''. The cardinality of a set $T$ is denoted as $\#(T)$. If $F$ is a field, then $\overline{F}$ doesn't denote the algebraic closure. For a morphisms of schemes $f:X\to S$,   the schematic preimage of $s\in S$ is $X_s$. Sheaves are denoted with the ``mathscr'' late$\chi$ font; particular the structure sheaf of a scheme $X$ is $\mathscr O_X$ (note the difference with the font $\mathcal O$). For any $\mathscr O_X$-module $\mathscr F$ and any $D\in\Div(X) $ we put $\mathscr F(D):=\mathscr F\otimes_{\mathscr O_X}\mathscr O_X(D)$. The notation $\det(\cdot)$ is used for the notion of ``determinant'' in the category of free modules over a ring and in the category of free $\mathscr O_X$-modules; the  exact meaning will be clear from the context.  If $K$ is a number field and $X\to\spec O_K$ is an integral scheme over the ring of integers $O_K$, then the function field of $X$ is denoted by $K(X)$. Finally it is important to point out that the letter $K$ will denote different mathematical objects in this paper (and in different contexts), so the reader should check at the beginning of each section its specific meaning from time to time.

\paragraph{Acknowledgements.} I would like to thank \emph{Ivan Fesenko} and \emph{Weronika Czerniaswka} from the university of Nottingham for the interesting discussion they had with me about the topic and their support. Moreover a special thanks to \emph{Fedor Bogomolov} and \emph{Nikolaos Diamantis} for reading my work and to  \emph{Dongwen Liu} for answering to my questions.

This work is supported by the EPSRC programme grant EP/M024830/1 (Symmetries and correspondences: intra-disciplinary developments and applications).

\section{Review of $2$-dimensional adelic geometry}\label{sect1}
\subsection{Abstract local theory}\label{subs1.1}
Let's recall the definition of local field:
\begin{definition}\label{locf}
A \emph{local field} (or a \emph{1-dimensional local field}) $F$, is one of the fields listed below:
\begin{itemize}
\item[$(1)$]$F=\mathbb R$ endowed with  the usual real absolute value$|\cdot|$.
\item[$(2)$] $F=\mathbb C$ endowed with  the usual complex  absolute value $||\cdot||$.
\item[$(3)$] $F$ is a complete discrete valuation field (the valuation is surjective) such that the residue field $\overline F$ is a perfect field. The valuation ring of $F$ is denoted as $\mathcal O_F$ and its maximal ideal is $\mathfrak p_F$. Moreover if $v$ is the valuation on $F$, then the absolute value  is given by $|x|_v:=q^{-v(x)}$, where $q=\#(\overline F)$ if $\overline F$ is a finite field, and $q=e:=\exp(1)$ otherwise.
\end{itemize}
If $F$ is of type $(1)$ or $(2)$, it is an \emph{archimedean local field} otherwise it is a \emph{non-archimedean local field}. A local field is topologized with the topology induced by the absolute value. A morphism between local fields is a continuous field homomorphism.
\end{definition}
\begin{remark}
According to our definition, a non-archimedean local field endowed with its natural topology is in general not locally compact.
\end{remark}
Remember that if $F$ is a non-archimedean local field there exists only one surjective complete valuation on it (see \cite[Theorem 1.4]{MM1}).  A higher local field is a simple generalization of definition $\ref{locf}$: given a complete discrete valuation field $F$, it might happen that the residue field $F^{(1)}:=\overline F$ is again a complete discrete valuation field; by taking one more time the residue field we have the field $F^{(2)}$. In other words, a complete discrete valuation field might originate a potentially infinite sequence of fields  $\{F^{(i)}\}_{i\ge 0}$ such that $F^{(0)}=F$ and $F^{(i+1)}=\overline {F^{(i)}}$. Each $F^{(i)}$ is called the \emph{$i$-th residue field}.
 
\begin{definition}
A \emph{$n$-dimensional local field}, for $n\ge 2$, is a complete discrete valuation field $F$ admitting sequence of residue fields $\{F^{(i)}\}_{i>0}$  such that  $F^{(n-1)}$ is a local field. If $F^{(n-1)}$ is an archimedean local field, then $F$ is called \emph{archimedean}, otherwise we say that $F$ is \emph{non-archimedean}. $F$ has \emph{mixed characteristic} if $\cha(F)\neq \cha({\overline F})$.
\end{definition}
\begin{example}\label{exlocf1}
The simplest $n$-dimensional local field is the field of iterated Laurent series over a perfect field $K$:
$$F=K((t_1))\ldots ((t_n))\,.$$
 If $f=\sum a_jt_n^j\in F$, with $a_i\in K((t_1))\ldots((t_{n-1}))$, we have the complete discrete valuation defined by $v(f)=\min\{j\colon a_j\neq 0\}$. The valuation ring is $\mathcal O_F=K((t_1))\ldots((t_{n-1}))[[t_n]]$ and the residue field is $F^{(1)}= K((t_1))\ldots((t_{n-1}))$. Clearly $F^{(n)}=K$.  When $n=2$, the elements of $K((t_1))((t_2))$ are the formal power series $\sum_{i,j} a_{i,j}t^i_1t_2^j$ such that $a_{i,j}=0$ when the indexes $i$ and $j$ are chosen in the following way: let's plot the couples $(j,i)$ as a lattice on the plane, then we select a semiplane like the one which is not coloured in figure \ref{lattice_coeff}. The coordinate $j$ is bounded from right, whereas the coordinate $i$ is bounded from above by a descending staircase line. 
 \begin{figure}[h]
 \centering
\begin{tikzpicture}
\draw[help lines, color=gray!70, dashed] (-4.9,-4.9) grid (4.9,4.9);
\draw[->,ultra thick] (-5,0)--(5,0) node[right]{$j$};
\draw[->,ultra thick] (0,-5)--(0,5) node[above]{$i$};

\def\myPath{(-4,5)
            --(-4,-1)
            --(-2,-1)
            --(-2,-3)
            --( 2,-3)
            --( 2,-4)
            --( 5,-4)}
    \draw \myPath;

    \begin{scope}[on background layer]
        \fill[gray!20] \myPath |- cycle;
    \end{scope}

\draw (-4,-4) node[right]{$a_{i,j}=0$};

\end{tikzpicture}
\caption{A cartesian diagram showing the lattice of couples $(j,i)$ corresponding to the coefficients $a_{i,j}$ of a power series in  $\sum_{i,j} a_{i,j}t^i_1t_2^j\in K((t_1))((t_2))$.}
\label{lattice_coeff}
\end{figure}
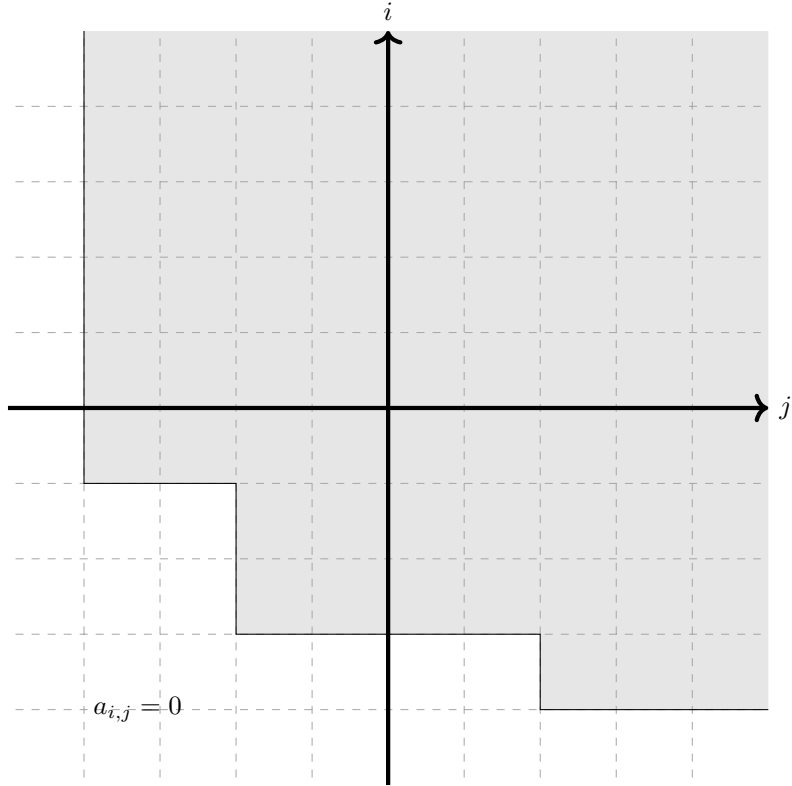
\end{example}

\begin{remark}
The above definition of dimension for a $n$-dimensional local field might seem quite counter-intuitive, indeed a $n$-dimensional local field can also be a $m$-dimensional local field for $m\neq n$. For instance $F=K((t_1))\ldots ((t_n))$  is $m$-dimensional for any $m=1,\ldots, n$. For our purposes it will be clear from the context which dimension we want to take in account. Often it is convenient to consider the maximum amongst all possible dimensions  (when it exists).
\end{remark}

\begin{remark}
Note in the case of archimedean $n$-dimensional local fields the $n$-th residue field doesn't exist.
\end{remark}

Let's give a less trivial example of higher local field:
\begin{example}\label{exlocf2}
Let $(K,v_K)$ be a non-archimedean local field and consider the following set of (double) formal series:
$$K\{\!\{t\}\!\}:=\left\{\sum_{j=-\infty}^{\infty} a_jt^j\colon a_j\in K,\; \inf_j v_K(a_j)>-\infty,\;\lim_{j\to-\infty}a_j=0\right\}$$
Addition and multiplication in $K\{\{t\}\}$ are defined in the following way:
\begin{equation}\label{sumhlf}
\sum_{j=-\infty}^{\infty} a_jt^j+\sum_{j=-\infty}^{\infty} b_jt^j=\sum_{j=-\infty}^{\infty} (a_j+b_j)t^j
\end{equation}
\begin{equation}\label{multhlf}
\sum_{j=-\infty}^{\infty} a_jt^j\cdot\sum_{j=-\infty}^{\infty} b_jt^j=\sum_{j=-\infty}^{\infty}\left(\sum_{r=-\infty}^{\infty} a_rb_{j-r}\right)t^j
\end{equation}
and $K\{\!\{t\}\!\}$  becomes a field. Note that the series with index $r$ in equation (\ref{multhlf}) is actually a convergent series in $K$. We can also define the following discrete valuation $v$ on $K\{\!\{t\}\!\}$:
\begin{equation}\label{valhlf}
v\left(\sum_{j=-\infty}^{\infty} a_jt^j\right):=\inf_j v_K(a_j)\,.
\end{equation}
It is not difficult to verify that $v$ is a well defined valuation and  $K\{\!\{t\}\!\}$ is complete with respect to $v$. Let's now analyse the structure of $F=K\{\!\{t\}\!\}$ as valuation field:
$$\mathcal O_F=\left\{\sum_{j=-\infty}^{\infty} a_jt^j\in K\{\!\{t\}\!\}\colon a_j\in \mathcal O_K\right\}$$
$$\mathfrak p_F=\left\{\sum_{j=-\infty}^{\infty} a_jt^j\in K\{\!\{t\}\!\}\colon a_j\in \mathfrak p_K\right\}$$
Consider the surjective homomorphism:
\begin{eqnarray*}
\pi:\;\; \mathcal O_F&\to&\overline{K}((t))\\
\sum a_jt^j &\mapsto& \sum\overline{a_j} t_j
\end{eqnarray*}
where clearly $\overline{a_j}$ is the natural image of $a_j$ in $\overline K$. Now it is evident that $\pi$ induces an isomorphism $\overline F\cong \overline{K}((t))$. In other words $F$ has a structure of $2$-dimensional local field such that $F^{(1)}=\overline{K}((t))$ and $F^{(2)}=\overline{K}$.  Clearly such a construction can be iterated several times to get the field:
$$K\{\!\{t_1\}\!\}\ldots\{\!\{t_n\}\!\}\,.$$ 
For example if $K=\mathbb Q_p$, then $K\{\!\{t\}\!\}$ is a 2-dimensional local field of mixed characteristic.
\end{example}
Remember that we have the following classical classification theorem for local fields:
\begin{theorem}[Classification theorem for local fields]\label{class1}
Let $F$ be a local field:
\begin{itemize}
\item[$(1)$] When $F$ is archimedean, then $F=\mathbb R$ or $F=\mathbb C$.
\item[$(2)$] When $F$ is not archimedean there are two cases:
\begin{itemize}
\item[$(2a)$] If $\cha F=\cha \overline F$, then $F\cong \overline F((t))$.
\item[$(2b)$] If $\cha F\neq \cha \overline F=p$, then $F$ is isomorphic to $K_p$ which denotes a finite extension of $ \mathbb Q_p$.
\end{itemize}
\end{itemize}
\end{theorem}
\proof
$(1)$ is true just by definition. For $(2)$ see for example \cite[II.5]{fesbook}.
\endproof
Such a classification can be extended for higher local fields, in particular any $n$-dimensional local field can be obtained by ``combining'' the higher local fields presented in examples \ref{exlocf1} and \ref{exlocf2}:
\begin{theorem}[Classification theorem for $n$-dimensional local fields]\label{classhigh}
Let $F$ be a $n$-dimensional local field with $n\ge 2$.
\begin{enumerate}
\item[$(1)$] If $\cha F=\cha F^{(1)}=\ldots=\cha F^{(n-1)}$, then 
$$F\cong F^{(n-1)}((t_1))\ldots((t_{n-1}))$$
and $F^{(n-1)}$ is isomorphic to one of the four fields listed in theorem \ref{class1}.
\item[$(2)$] If $r\in\{2,3,\ldots n\}$ is  the unique number such that $\cha F^{(n-r)}\neq\cha F^{(n-r+1)}=p$, then: 
\begin{enumerate}  
\item[$(2a)$] When $r\neq n$,  $F$ is isomorphic to a finite extension of:
$$
K_p\{\!\{t_1\}\!\}\ldots \{\!\{t_{r-1}\}\!\}((t_{r}))\ldots((t_{n-1}))\,.
$$
\item[$(2b)$] When $r=n$ (i.e. in the mixed characteristic case),  $F$ is isomorphic to a finite extension of:
$$
K_p\{\!\{t_1\}\!\}\ldots \{\!\{t_{n-1}\}\!\}\,.
$$
\end{enumerate}
\end{enumerate}
\end{theorem}
\proof
See \cite[Theorem 2.18]{MM1}.
\endproof
In this paper we will focus mainly on $2$-dimensional local fields, so let's give  a table with all possible $2$-dimensional local fields by using the classification theorem:

\begin{equation}\label{2dimtable}
\begin{centering}
\begin{tabular}{ |p{2.7cm}|p{2.7cm}|p{2.7cm}|p{2.7cm}|}
 \hline
 \multicolumn{4}{|c|}{$2$-dimensional local fields} \\
 \hline
   \multicolumn{1}{|c|}{Geometric} & \multicolumn{2}{|c|}{Arithmetic} &\multicolumn{1}{|c|}{Archimedean}\\
\hline
 $(0,0,0)$, $(p,p,p)$  & $(0,p,p)$ & $(0,0,p)$ & \\
\cline{1-3}
 &&&\\
 $K((t_1))((t_2))$ with $K$ perfect  & finite extension of $K_p\{\!\{t\}\!\}$&   $K_p((t))$& $\mathbb C((t))$ or $\mathbb R((t))$\\
 &&&\\
 \hline
\end{tabular}
\end{centering}
\end{equation}
\\

For a non-archimedean local field $(F,v)$, we have the notion of \emph{local parameter} $\varpi$ which is any generator of the maximal ideal $\mathfrak p_F$, or equivalently any element such that $v(\varpi)=1$. Clearly we have the (recursive) generalization for $n$-dimensional local fields.
\begin{definition}
Let $F$ be a non-archimedean $n$-dimensional local field, then a \emph{sequence of local parameters for $F$} is a $n$-tuple $(\varpi_1,\ldots,\varpi_n)\in  F\times\ldots\times F$ satisfying the following properties:
\begin{itemize}
\item[\sbt] $\varpi_n$ is a local parameter for $F$.
\item[\sbt] $(\varpi_1,\ldots,\varpi_{n-1})\in\mathcal O_F\times\ldots\times \mathcal O_F$ and the sequence of natural projections $(\overline{\varpi}_1,\ldots,\overline{\varpi}_{n-1})$ is a sequence of local parameters for the residue field $\overline F$.
\end{itemize}
\end{definition}
One can obtain a sequence of local parameters, by applying the following algorithm: choose any local parameter for $F^{(n-1)}$, then pick any of its liftings in $F$, this will be $\varpi_1$.  Choose choose any local parameter for $F^{(n-2)}$, then pick any of its liftings in $F$, this will be $\varpi_2$, etc. Let's give another basic  definition: 
\begin{definition}
Let $F$ be a $n$-dimensional local field and put  $\mathcal O^{(0)}_F:=F$, then we define recursively the \emph{$j$-th valuation ring} (for $j\ge1$):
$$\mathcal O^{(j)}_F:=\left\{x\in\mathcal O_F\colon \overline{x}\in\mathcal O^{(j-1)}_{\overline F}\right\}$$
It is clear that $\mathcal O^{(1)}_F=\mathcal O_F$. For the algebraic properties of $\mathcal O^{(j)}_F$ the reader can check \cite[3]{MM1}.
\end{definition}

 From now on in this section we assume that $F$ is a $2$-dimensional local field such that:
\begin{enumerate}
\item[\sbt] $\cha F=0$ and $\cha F^{(2)}=p$.
\item[\sbt] $F$ is endowed with a ST-ring\footnote{A ST-ring (semi-topological ring) is just a ring endowed with a linear topology as additive group such that the multiplication for any fixed element is continuous. Note that a ST-ring is not necessarily a topological ring. See \cite{yek1} for more details.} topology and there exists a mixed characteristic local field $K$ with a fixed embedding $K\hookrightarrow F$ of ST-rings ($K$ has the discrete valuation topology).
\end{enumerate}
In this case we say that $F$ is an \emph{arithmetic $2$-dimensional local field over $K$}. The presence of the local field  $K$ inside $F$ comes from the theory of arithmetic surfaces and it will be explained in section \ref{sect2}. 

\paragraph{Equal characteristic.}
Let $F$ be an arithmetic $2$-dimensional local field such that $\cha \overline F=0$. 
\begin{definition}
The \emph{coefficient field of $F$ (with respect to $K$)} is the algebraic closure of $K$ inside $F$ and it is denoted as $k_F$.
\end{definition}
The coefficient field $k_F$ is a finite extension of $K$ and moreover $\overline F=k_F$. In particular  $F\cong k_F((t))$. The valuation field $F$ is naturally endowed with the usual tame symbol $(\,,\,)_{F}: F^\times\times F^\times\to k_F^\times$, so we can obtain the Kato symbol (or two dimensional tame symbol) by simply composing it with the field norm map:
\begin{definition}
The \emph{Kato symbol} for $F$ (with respect to $K$) is given by:
$$(\,,\,)_{F|K}:N_{k_F|K}\circ (\,,\,)_{F}:F^\times\times F^\times\to K^\times\,.$$
\end{definition}

\paragraph{Mixed characteristic.} 
Now we assume that $F$ is an arithmetic $2$-dimensional local field of mixed characteristic. By the classification theorem $\mathbb Q_p$ is contained in $F$ and we have the notion of constant field of $F$ which replaces the one of coefficient field:
\begin{definition}
The \emph{constant field of $F$} is the algebraic closure of $\mathbb Q_p$ in $F$, and it is denoted by $k_F$.
\end{definition}
\begin{remark}
Note that the definition of the constant field doesn't depend on $K$ so it makes sense for any $2$-dimensional local field of mixed characteristic. Of course it might happen that $K=\mathbb Q_p$. 
\end{remark}
Since $K$ is a finite extension of $\mathbb Q_p$ (by the $1$-dimensional classification theorem), we know that $k_F$ is an intermediate field between $K$ and $F$. The constant field $k_F$ is a finite extension of $\mathbb Q_p$ (so also a finite extension of $K$).
\begin{definition}
We say that an arithmetic $2$-dimensional local field of mixed characteristic $F$ is \emph{standard} if there is a $k_F$ isomorpshism $F\cong k_F\{\!\{t\}\!\}$. When an isomorphism is given, we say that we have fixed a \emph{parametrization} of $F$.
\end{definition}
We will study standard fields first and extend any result for a generic $F$ thanks to the following result:
\begin{proposition}\label{stand}
There exists a standard field $L$ contained in $F$ such that: $[F:L]<\infty$, $k_F=k_L$ and $\overline F=\overline L$.
\end{proposition}
\proof
See \cite[Lemma 2.14]{MM3}.
\endproof
So, from now on in this subsection we fix $L$ to be a standard field contained in $F$ with the properties described in proposition \ref{stand}. Clearly we have the following field extensions that need to be kept always in mind (we mark the finite extensions with the superscript $\mathfrak f$):
\begin{equation}\label{field_struct}
\mathbb Q_p\subseteq^{\mathfrak f} K\subseteq^{\mathfrak f} k_L=k_F\subseteq L\cong k_L\{\!\{t\}\!\}\subseteq^{\mathfrak f} F\,.
\end{equation}

Finally, we want to define the Kato symbol for $F$ and the strategy is the usual one: we start from $k_L\{\!\{t\}\!\}$ and we extend our arguments to $F$ by checking that everything is independent from parametrizations and from the choice of the standard fields. We will heavily use some $K$-theoretic notions developed in appendix \ref{k_th}.

Fix just for the moment $L=k_L\{\!\{t\}\!\}$, then we define:
\begin{equation}
\begin{tikzcd}
(\,,\,)_{L|K}: L^\times\times L^\times\arrow[r, "{\{\,,\,\}}"] & K_2(L)\arrow[r] & \widehat{K}_2(L)\arrow[r, "-\res^{(2)}_L"] & \widehat{K}_1(k_L)=k_L^\times\arrow[r,"N_{k_L|K}" ] & K^\times\\
\end{tikzcd}
\end{equation}
where:
\begin{enumerate}
\item[\sbt] $\{\,,\,\}$ is the natural projection arising from the definition of $K_2(L)$ (see proposition \ref{K-gr}).
\item[\sbt] The morphism  $K_2(L)\to\widehat{K}_2(L)$ is the map given by the construction of $\widehat{K}_2(L)$ as projective limit (see equation (\ref{completedK})).
\item[\sbt] $\res^{(2)}_L$ is the higher Kato residue map constructed in theorem \ref{Kato_res}. Note that $\widehat{K}_1(k_L)=k_L^\times$ because $k_L$ is already complete.
\end{enumerate} 
Moreover by simplicity we use the following notation:
\begin{equation}\label{not_boundarymap}
\begin{tikzcd}
\partial_L: K_2(L)\arrow[r] & \widehat{K}_2(L)\arrow[r, "-\res^{(2)}_L"] & k_L^\times\,.\\
\end{tikzcd}
\end{equation}

\begin{remark}
\cite{otherLiu} gives an explicit description of $\res^{(2)}_L$ which involves winding numbers.
\end{remark} 
\begin{definition}
Let $L$ be a generic standard field and fix a parametrization: $p: L\to k_L\{\!\{t\}\!\}$  then we define:
$$
\begin{tikzcd}
(\,,\,)_{L|K}: L^\times\times L^\times\arrow[r, "{\{\,,\,\}}"] & K_2(L)\arrow[r, "K_2(p)"] & K_2(k_L\{\!\{t\}\!\})\arrow[r, "\partial_{k_L\{\!\{t\}\!\}}"] & k_L^\times\arrow[r,"N_{k_L|K}" ] & K^\times\\
\end{tikzcd}
$$
and we put $\partial_L:=\partial_{k_L\{\!\{t\}\!\}}\circ K_2(p)$.
\end{definition}
\begin{proposition}\label{not_para}
Let $L$ be a standard field, then the definition of $(\,,\,)_{L|K}$ doesn't depend on the parametrization of $L$.
\end{proposition}
\proof
See \cite[Corollary 3.7]{otherLiu}.
\endproof
At this point we are ready to give the general definition of the Kato symbol:
\begin{definition}
Let $F$ an arithmetic $2$-dimensional local field and let $L$ be a standard field contained in $F$, then the \emph{Kato symbol for $F$} (with respect to $K$) is given by:
\begin{equation}
\begin{tikzcd}
(\,,\,)_{F|K}: F^\times\times F^\times\arrow[r, "{\{\,,\,\}}"] & K_2(F)\arrow[r, "K_2(N_{F|L})"] & K_2(L)\arrow[r, "\partial_L"] & k_L^\times\arrow[r, "N_{k_L|K}"] & K^\times\\
\end{tikzcd}
\end{equation}
\end{definition}
\begin{proposition}
The definition of $(\,,\,)_{F|K}$ doesn't depend on the choice of $L$ inside $F$.
\end{proposition}
\proof
See \cite[Proposition 3]{kato}.
\endproof
\subsection{Adelic geometry}\label{sect1.2}
Let's fix $B=\spec O_K$ for a number field $K$;  $\varphi:X\to B$ is a $B$-scheme satisfying the following properties:
\begin{enumerate}
\item[\sbt] $X$ is  two dimensional, integral, and regular. The generic point of $X$ is $\eta$ and the function field of $X$ is denoted by $K(X)$.
\item[\sbt] $\varphi$ is proper and flat.
\item[\sbt] The generic fibre, denoted by $X_K$, is a geometrically integral, smooth, projective curve over $K$. The generic point of $B$ is denoted by $\xi$.
\end{enumerate} 
We say that $X$ is an \emph{arithmetic surface over $B$}. We consider the set of all possible \emph{flags} $x\in Y\subset X$ where $x$ is a closed point of $X$ contained in an integral curve $Y$.

From now on a curve $Y$ on $X$ will always be an integral curve and its unique generic point will be denoted with the letter $y$.  By simplicity we will often identify  $Y$ with its generic point $y$, which means that by an abuse of language  and notation we will use sentences like ``let $y\subset X$ be a curve on $X$..." or ``let $x\in y\subset X$ be a flag on $X$...''. In other words $y$  is considered  as a scheme or as a point depending on  the context.
\begin{definition}
Fix a closed point $x\in X$, then:
\begin{itemize}
\item[\sbt] $\mathcal O_x:=\widehat{\mathscr O_{X,x}}$. It is a Noetherian, complete, regular, local, domain of dimension $2$ with maximal ideal $\widehat{\mathfrak m_x}$.
\item[\sbt] $K'_x:=\fr\mathcal O_x$.
\item[\sbt] $K_x:=K(X)\mathcal O_x\subseteq K'_x$. Notice that this is not a field.
\end{itemize}
For a curve $y\subset X$ we put:
\begin{itemize}
\item[\sbt] $\mathcal O_y:=\widehat{\mathscr O_{X,y}}$. It is a complete DVR  with maximal ideal $\widehat{\mathfrak m_y}$.
\item[\sbt] $K_y:=\fr\mathcal O_y$. It is a complete discrete valuation field with valuation ring $\mathcal O_y$. The valuation is denoted by $v_y$.
\end{itemize}
For any point $b\in B$ we put:
\begin{itemize}
\item[\sbt] $\mathcal O_b:=\widehat{\mathscr O_{B,b}}$. It is a complete DVR.
\item[\sbt] $K_b:=\fr \mathcal O_b$. It is a local field with finite residue field. The valuation is denoted by $v_b$.
\end{itemize}
\end{definition}
Fix a flag $x\in y\subset X$, then we have a surjective local homomorphism $\mathscr O_{X,x}\to\mathscr O_{y,x}$ with kernel $\mathfrak p_{y,x}$ induced by the closed embedding $y\subset X$ (note that $\mathfrak p_{y,x}$ is a prime ideal of height $1$).

 The inclusion $\mathscr O_{X,x}\subset \mathcal O_x$ induces a morphism of schemes $\phi:\spec\mathcal O_x\to\spec\mathscr O_{X,x}$ and we define the \emph{local branches of $y$ at $x$} as the elements of the set
$$y(x):=\phi^{-1}(\mathfrak p_{y,x})=\set{\mathfrak z\in\spec\mathcal O_x\colon \mathfrak z\cap \mathscr O_{X,x}=\mathfrak p_{y,x}}\,.$$
If $y(x)$ contains only an element, we say that $y$ is unbranched at $x$. 

\begin{figure}[ht]
\centering
\begin{tikzpicture}
[spy using outlines={circle, magnification=4, size=2.9cm, connect spies}]

\draw[name path=c1]  (0,-1.5)  .. controls (0.5,0.5) and (1,0.5) .. (1,0);
\draw[name path=c2] (1,0) .. controls (1,-0.5) and (0.5,-0.5) .. (0,1.5)node[anchor=west]{$y$}; 
\draw (-3,2)..controls (-1,2.5) and (1,1.5)..(3,2)node[anchor=west]{$X$};
\draw (-3,2)--(-3,-2);
\draw (-3,-2)..controls (-1,-1.5) and (1,-2.5)..(3,-2);
\draw (3,2)--(3,-2);

\fill [name intersections={
of=c1 and c2,  by={a,b}}]
(a)
(b) circle (1.5pt) node[anchor=east]{$x$};

\spy[lightgray] on (b) in node  at (1,-4.5);
\draw (2.9,-4.8)node{$y(x)$};
\end{tikzpicture}
\caption{\footnotesize{Informally the local branches of $y$ at $x$ can be depicted in the following way: consider a small neighbourhood of $x$, then each distinct ``piece of $y$'' that we see passing through $x$ corresponds to a local branch $\mathfrak z$. In this particular case $y$ has a simple node at $x$, so $2$ local branches at $x$.}}
\end{figure}
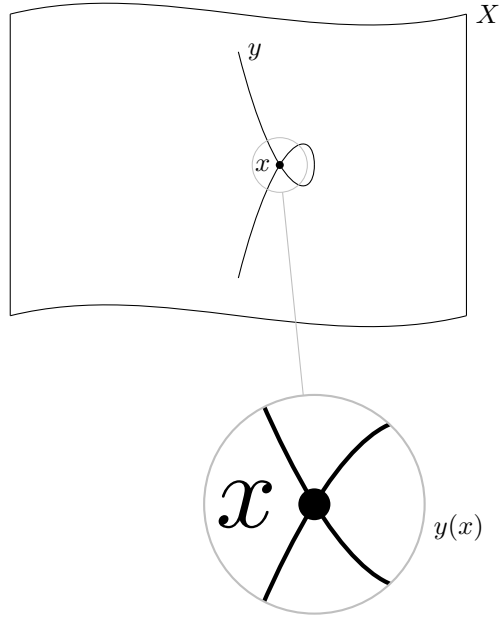
\begin{remark}
If $x$ is a cusp point on $y$, one can show that $y$ unbranched at $x$.
\end{remark}

\begin{definition}
Let $\mathfrak z\in y(x)$ be a local branch of a curve $y$ at point $x$, then let's define the field
$$K_{x,\mathfrak z}:=\fr\left(\widehat{\left(\mathcal O_x\right)_\mathfrak z}\right)\,.$$
in other words: we localise $\mathcal O_x$ at the prime ideal $\mathfrak z$, then we complete it at its maximal ideal and finally we take the fraction field. By convenience we put $\mathcal O_{x,\mathfrak z}:=\widehat{\left(\mathcal O_x\right)_\mathfrak z}$.
\end{definition}

The proof of the following proposition relies on some basic commutative algebra results:
\begin{proposition}\label{two_dim_loc}
Let $x\in y\subset X$ be a flag and let $\mathfrak z\in y(x)$. Then  $K_{x,\mathfrak z}$ is a $2$-dimensional valuation field  such that $\mathcal O_{K_{x,\mathfrak z}}=\mathcal O_{x,\mathfrak z}$ and $K^{(2)}_{x,\mathfrak z}$ is a finite extension of $k(x)$.
\end{proposition}
\proof
First of all $\hei \mathfrak z\ge \hei \mathfrak p_{y,x}=1$, but if $\hei \mathfrak z=2$ then $\mathfrak z$ is the maximal ideal of $\mathcal O_x$ and we have that $\mathfrak z\cap
\mathscr O_{X,x}=\mathfrak m_x$, a contradiction. Therefore $\hei \mathfrak z=1$ and $\dim\left(\mathcal O_x\right)_{\mathfrak z}=1$. It follows that $\widehat{\left(\mathcal O_x\right)_\mathfrak z}$ is a Noetherian, complete, local, domain of dimension $1$, i.e. a complete DVR which is the valuation ring of the complete discrete valuation field $K_{x,\mathfrak z}$. The residue field of $K_{x,\mathfrak z}$ is by definition:
$$K^{(1)}_{x,\mathfrak z}:=\left(\mathcal O_x\right)_\mathfrak z/\mathfrak z\left(\mathcal O_x\right)_\mathfrak z=\fr\left(\mathcal O_x/\mathfrak z\right)\,.$$
Note that $\mathcal O_x/\mathfrak z$ is a Noetherian, complete, local domain of dimension $1$ (in general we may lose the regularity by passing to the quotient). Consider the normalisation $\widetilde {\mathcal O_x/\mathfrak z}$ of $\mathcal O_x/\mathfrak z$; the domain $\widetilde {\mathcal O_x/\mathfrak z}$ is obviously normal and again Noetherian and complete. Moreover by Nagata theorem (see \cite[Ch. IX, 4, no 2, Theorem 2]{bou}) $\mathcal O_x/\mathfrak z$ is a Japanese ring, therefore in particular  $\widetilde {\mathcal O_x/\mathfrak z}$ is a finite $\mathcal O_x/\mathfrak z$-module. Now \cite[Corollary 7.6]{Ei} implies that $\widetilde {\mathcal O_x/\mathfrak z}$ is also local, and by summing up all the listed property we can conclude that $\widetilde {\mathcal O_x/\mathfrak z}$ is a complete DVR with fraction field $\fr\left(\mathcal O_x/\mathfrak z\right)$. This proves that $K^{(1)}_{x,\mathfrak z}$ is a complete valuation field.\\
It remains to show only that the second residue field $K^{(2)}_{x,\mathfrak z}$ is a finite extension of $k(x)$. By definition $K^{(2)}_{x,\mathfrak z}$ is the residue field of the local ring $\widetilde {\mathcal O_x/\mathfrak z}$, but we already know that $\widetilde {\mathcal O_x/\mathfrak z}$ is a finite $\mathcal O_x/\mathfrak z$-module, so  $K^{(2)}_{x,\mathfrak z}$ is a finite extension of:

$$\left(\mathcal O_x/\mathfrak z\right)\big/\left(\widehat{\mathfrak m_x}/\mathfrak z \right)\cong \mathcal O_x/\widehat{\mathfrak m_x}\cong\mathscr O_{X,x}/\mathfrak m_x=k(x)\,.$$
\endproof

It is not amongst the purposes of this paper to treat the topology of $K_{x,\mathfrak z}$, but it is enough to know that there are several ways to topologise $K_{x,\mathfrak z}$, some of them are equivalent, and we end up with  a structure of ST-ring on $K_{x,\mathfrak z}$. See for example \cite[1.]{Brau} for a survey about topologies on $K_{x,\mathfrak z}$.

\begin{definition}
Let $x\in y\subset X$ be a flag and let $\mathfrak z\in y(x)$, then we put $E_{x,\mathfrak z}:=K^{(1)}_{x,\mathfrak z}$ and $k_\mathfrak z (x):=K^{(2)}_{x,\mathfrak z}$. 
$$
\begin{tikzcd}
K_{x,\mathfrak z}\arrow[r, phantom, "\supset"]\arrow[dr, bend right,dashed, no head
] & \mathcal O_{x,\mathfrak z}:=\mathcal O_{K_{x,\mathfrak z}}\arrow[d]\arrow[r, phantom, "\supset"]& \mathcal O^{(2)}_{x,\mathfrak z}:=\mathcal O_{K_{x,\mathfrak z}}^{(2)}\arrow[d]\\
& E_{x,\mathfrak z}:=K^{(1)}_{x,\mathfrak z}\arrow[r, phantom, "\supset"]\arrow[dr, bend right,dashed, no head
] & \mathcal O_{E_{x,\mathfrak z}}\arrow[d]\\
& & k_\mathfrak z(x):=K^{(2)}_{x,\mathfrak z}
\end{tikzcd}
$$
The valuation on $K_{x,\mathfrak z}$ is $v_{x,\mathfrak z}$ and the valuation on $E_{x,\mathfrak z}$ is $v^{(1)}_{x,\mathfrak z}$.
Moreover:
$$K_{x,y}:=\prod_{\mathfrak z\in y(x)}K_{x,\mathfrak z}\,,\quad\mathcal O_{x,y}:=\prod_{\mathfrak z\in y(x)}\mathcal O_{x,\mathfrak z}\,,\quad\mathcal O^{(2)}_{x,y}:=\prod_{\mathfrak z\in y(x)}\mathcal O^{(2)}_{x,\mathfrak z}\,,$$
$$E_{x,y}:=\prod_{\mathfrak z\in y(x)}E_{x,\mathfrak z}\,,\quad k_y(x):=\prod_{\mathfrak z\in y(x)}k_\mathfrak z (x)\,.$$
\end{definition}
\begin{remark}\label{commremark}
Remember from commutative algebra  the following chain of implications:
$$(\text{$A$ regular local})\Rightarrow (\text{$A$ a UFD})\Rightarrow(\text{Any prime $\mathfrak p$ s.t. $\hei(\mathfrak p)=1$ is principal})\,.
$$
So, $\mathscr O_{X,x}$ is a UFD and $\mathfrak p_{y,x}$ is principal, but also $\mathcal O_x$ is a UFD and $\mathfrak z$ is principal.
\end{remark}
\begin{proposition}\label{localpar}
Let $\mathfrak p_{y,x}=(\varpi_y)$ for $\varpi_y\in\mathscr O_{X,x}$, then we can choose the uniformizing parameter for $K_{x,\mathfrak z}$ to be $\varpi_y$. 
\end{proposition}
\proof We show that $\varpi_y$ generates the maximal ideal of $\mathcal O_{x,\mathfrak z}$. First of all we notice that the ring $\mathscr O_{X,x}/\varpi_y\mathscr O_{X,x}\cong\mathscr O_{y,x}$ is reduced, and this implies that $\widehat{\mathscr O_{y,x}}=\mathcal O_x/\varpi_y\mathcal O_x$ is reduced too. By remark \ref{commremark} $\varpi_y$ has a unique factorization $\varpi_y=p_1\ldots p_m$ in $\mathcal O_x$, and all the $p_i's$ are distinct prime elements thanks to the fact that $\mathcal O_x/\varpi_y\mathcal O_x$ is reduced. Again remark \ref{commremark} implies that $\mathfrak z=(p_j)$ for some index $j$. Any element of $\mathfrak z(\mathcal O_x)_{\mathfrak z}$ can be written as $\frac{p_ja}{b}$ with $b\notin\mathfrak z$ but:
$$\frac{p_ja}{b}=\frac{p_1\ldots p_m a}{p_1\ldots p_{j-1}p_{j+1}\ldots p_m b }=\frac{\varpi_ya}{p_1\ldots p_{j-1}p_{j+1}\ldots p_m b}$$
Since $p_1\ldots p_{j-1}p_{j+1}\ldots p_m b\notin\mathfrak z$, we can conclude that $\varpi_y$ generates the prime ideal $\mathfrak z(\mathcal O_x)_{\mathfrak z}$ of $(\mathcal O_x)_{\mathfrak z}$.
\endproof
\begin{corollary}\label{choosepam}
If $\varpi_y$ is a uniformizing parameter for the complete valuation field $K_y$, then it is a uniformizing parameter for $K_{x,\mathfrak z}$.
\end{corollary}
\proof
 It follows from  proposition \ref{localpar} and the fact that $(\mathscr O_{X,x})_{\mathfrak p_{y,x}}\cong\mathscr O_{X,y}$.
\endproof

\begin{remark}\label{first_res_rem}
Fix a flag  $x\in y\subset X$. The local homomorphism $\mathscr O_{X,x}\to\mathscr O_{y,x}$ induces a local homomorphism $\mathcal O_x\to\widehat{\mathscr O_{y,x}}$ which gives a bijective correspondence between the ideals in $y(x)$ and the minimal prime ideals of $\widehat{\mathscr O_{y,x}}$. Hence the ring of adeles of the curve $y$ is recovered in the following way:
$$\mathbf A_y=\sideset{}{'}\prod_{x\in y} E_{x,y}\,.$$
\end{remark}
\begin{proposition}\label{restrval}
Let's denote with $v_{x,\mathfrak z}$ the valuation of $K_{x,\mathfrak z}$ and with $v_y$ the valuation of $K_y$. Then the restriction of $v_{x,\mathfrak z}$ to $K_y$ is equal to $v_y$.
\end{proposition}
\proof
By remark \ref{first_res_rem} we deduce that $E_{x,\mathfrak z}$ contains $k(y)$, which is in turns the residue field of $K_y$, so the claims follows directly from corollary \ref{choosepam}.
\endproof

 The structure of $K_{x,\mathfrak z}$ depends on the nature of the curve $y$ and we can distinguish two cases:

\paragraph{$y$ is a vertical curve.} If $\varphi(y)=b\in B$, then $y$ is a projective curve over the finite field $k(b)$; we assume that $k(b)$ has characteristic $p$.  $K_{x,\mathfrak z}$ has characteristic $0$ since we have the embeddings $\mathbb Q\subset K\subset K(X)\subset K_{x,\mathfrak z}$ and the residue field $E_{x,\mathfrak z}$ has characteristic $p$ since $k(b)\subset k(y)\subset E_{x,\mathfrak z}$. We conclude that $K_{x,\mathfrak z}$ is a two dimensional local field of type $(0,p,p)$ and by the classification theorem we have that $K_{x,\mathfrak z}$ is a finite extension of $K_p\{\!\{t\}\!\}$ where  $K_p$ is a finite extension of $\mathbb Q_p$.

\paragraph{$y$ is a horizontal curve.} In this case $K_{x,\mathfrak z}$ has still characteristic $0$, but we have the embedding $K\subseteq k(y)$ given by the surjective map $y\to B$; therefore $E_{x,\mathfrak z}$ has characteristic $0$. Moreover, if $\varphi(x)=b$, the local homomorphism $\varphi^{\#}_x:\mathscr O_{B,b}\to\mathscr O_{X,x}$ induces a field embedding $k(b)\subseteq k(x)$ and this implies that $k_{\mathfrak z}(x)$ has characteristic $p$. We conclude that $K_{x,\mathfrak z}$ is a two dimensional local field of type $(0,0,p)$ and by the classification theorem we have that $K_{x,\mathfrak z}\cong K_p((t))$.  \\

If $\varphi(x)=b$ we have an induced embedding $K_b\hookrightarrow K_{x,\mathfrak z}$, so we can conclude that $K_{x,\mathfrak z}$ is an arithmetic $2$-dimensional local field over $K_b$ and we can apply the local theory developed in subsection \ref{subs1.1}. 

The ring of adeles $\mathbf A_X$ will be the result of a ``glueing'' of the local data $\{K_{x,y}\}_{x\in y\subset X}$ where the couple $(x,y)$ runs amongst all flags in $X$. The glueing procedure  will be described precisely, but roughly speaking we will define  $\mathbf A_X$ inside  the big product of rings
$$\mathbf A_X\subset \prod_{\substack {x\in y,\\ y\subset X}}K_{x,y}$$
as a sort of ``double restricted product''.\\

\paragraph{First ``restricted product'': the  rings $\mathbb A_y^{(r)}$ and  $\mathbb A_y$.}
In this paragraph we fix a curve $y\subset X$, and denote with $\mathfrak J_{x,y}$ the Jacobson radical of $\mathcal O_{x,y}$. 
\begin{definition}
Let's put:
\begin{equation*}
\mathbb A^{(0)}_y = \left\{ 
  \begin{aligned}
  & (\alpha_{x,y})_{x\in y}\in\prod_{x\in y}\mathcal O_{x,y} \colon  \forall s >0,\; \alpha_{x,y}\in\mathcal O_x+\mathfrak J_{x,y}^s \\ 
  & \text{for all but finitely many $x\in y$.}
  \end{aligned}
\right\}\subset \prod_{x\in y} \mathcal O_{x,y}
\end{equation*}
then for any $r\in\mathbb Z$
$$\mathbb A^{(r)}_y:= \widehat{\mathfrak m}^{r}_{y}\mathbb A_y^{(0)}\subset \prod_{x\in y} K_{x,y}\,.$$
Clearly  $\mathbb A^{(r)}_y\supseteq\mathbb A^{(r+1)}_y$ and $\bigcap_{r\in\mathbb Z} \mathbb A^{(r)}_y=0$. Moreover we  define 
$$\mathbb A_y:=\bigcup_{r\in\mathbb Z} \mathbb A^{(r)}_y\,.$$
\end{definition}
\begin{remark}
We have the inclusion $\mathbb A_y\subset\prod_{x\in y} K_{x,y}$, therefore we can interpret  $\mathbb A_y$ as a ``restricted product'' of the rings $K_{x,y}$ for $y$ fixed and $x\in y$. Thus we can write:
 $$\mathbb A_y=\sideset{}{'}\prod_{x\in y} K_{x,y}$$
 where $\sideset{}{'}\prod$ here is just a piece of notation without any formal meaning.
 \end{remark}

Each $\mathbb A_y^{(r)}$ can be endowed with a  ind/pro linear topology, and $\mathbb A_y$ can be seen as linear direct limit of the topological groups $\mathbb A_y^{(r)}$. More details about such topologies will be given in the second paper of the series.

\paragraph{Second ``restricted product'': the ring $\mathbf A_X$.}
The construction of $\mathbb A_y$ can be seen as a  way to take the restricted product of $\prod_{x\in y} K_{x,y}$. The final step in order to construct the ring of adeles $\mathbf A_X$ is to take the restricted product of the groups $\mathbb A_y$ over all the curves in $X$ with respect to the subgroups $\mathbb A_y^{(0)}$.
\begin{definition}
$$\mathbf A_X:=\left\{(\beta_y)_{y\subset X}\in\prod_{y\subset X}\mathbb A_y\colon \beta_y\in\mathbb A^{(0)}_y\; \text{for all but finitely many $y$}\right\}\subset\prod_{\substack {x\in y,\\ y\subset X}} K_{x,y}\,.$$
In a more suggestive way, we write by commodity
$$\mathbf A_X=\sideset{}{''}\prod_{\substack {x\in y\\ y\subset X}} K_{x,y}$$
where the symbol ``$\sideset{}{''}\prod$'' is just a piece of notation which remembers  that we are taking a ``double restricted product''.
\end{definition}

\begin{remark}
It is fundamental to recall that $\mathbf A_X$ is \emph{not} the full ring of adeles associated to the completed surface $\widehat X$, bcause we didn't take in account the fibres at infinity.
\end{remark}
In order to topologise  $\mathbf A_X$ we need to recall the description of the restricted product, by means of categorical limits, for linearly topologised groups. Let $\{G_i\}_{i\in I}$ a set of linearly topologised groups and for any $i$ let $H_i\subset G_i$ be a subgroup endowed with the subspace topology. We denote the family of finite subsets of $I$ as $\mathcal P_f(I)$; it forms a directed set with the relation $J\subseteq J'$. For any $J\in\mathcal P_f(I)$ define
$$G_J:=\prod_{i\in J} G_i\times\prod_{i\not\in J} H_i\,,$$
if  $J\subseteq J'$  the identity in each factor induces an embedding $G_J\hookrightarrow G_{J'}$, thus we have a direct system $\{G_J\}_J$ and it is easy to see that 
$$\sideset{}{'}\prod_i G_i=\varinjlim_J G_J\,,$$
where $\sideset{}{'}\prod_i G_i$ is the usual restricted product of the $G_i$ with respect to the subgroups $H_i$. At this point, on each $G_J$ we put the product topology and $\sideset{}{'}\prod_i G_i$ is endowed with  the linear  direct limit topology. By definition $\mathbf A_X$ is the restricted product of the groups $\mathbb A_y$ with respect to the subgroups  $\mathbb A_y^{(0)}$ for any $y\subset X$. Therefore we endow  $\mathbf A_X$ with the topology described above.

 We now introduce some important subspaces in order to construct the  adelic complexes associated to the surface $X$. Here the definitions are made ``by hands'', but such subspaces can be recovered as a particular case of the general theory of Beilinson adeles (see \cite[8]{MM1}). First of all let's consider the following diagonal embeddings:
 $$K_x\subset\prod_{y\ni x} K_{x,y}, \quad K_y\subset\prod_{x\in y} K_{x,y}\,,$$
so  we can consider:
$$\prod_{x\in X}K_x\subset\prod_{\substack {x\in y\\ y\subset X}} K_{x,y}, \quad \prod_{y\subset X}K_y\subset\prod_{\substack {x\in y\\ y\subset X}} K_{x,y}\,.$$
Let's define:
$$A_{012}:=\mathbf A_X\,;\quad A_{12}:=\mathbf A_X\cap \prod_{\substack{x\in y\\y\subset X}}\mathcal O_{x,y}=\prod_{y\subset X} \mathbb A^{(0)}_y; $$
$$A_{02}:=\mathbf A_X\cap\prod_{x\in X} K_x\,;\quad A_{2}:=\mathbf A_X\cap\prod_{x\in X} \mathcal O_x\,;\quad A_{01}:=\mathbf A_X\cap\prod_{y\subset X} K_y\,;$$
$$\quad A_{1}:=\mathbf A_X\cap\prod_{y\subset X} \mathcal O_y\,;\quad A_0:=K(X)$$
The containment relations are depicted in the following diagram:
$$
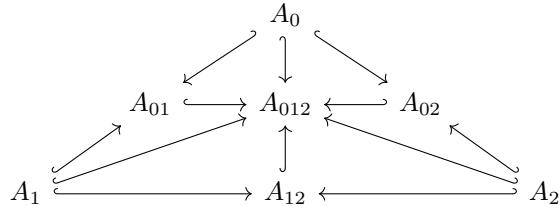
\begin{tikzcd}
&&A_0\arrow[d,hook']\arrow[dl,hook']\arrow[dr,hook]&&\\
&A_{01}\arrow[r,hook]& A_{012} &A_{02}\arrow[l,hook']&\\
A_1\arrow[rr,hook]\arrow[ur,hook]\arrow[urr,hook]&&A_{12}\arrow[u,hook]&&A_2\arrow[ll,hook']\arrow[ul,hook']\arrow[ull,hook']\\
\end{tikzcd}
$$
and we have the adelic complex:
\begin{equation}\label{ad_coomplex2}
\begin{tikzcd}
\mathcal A_X: & A_0\oplus A_1\oplus A_2 \arrow[r, "d^0"]& A_{01}\oplus A_{02}\oplus A_{12}\arrow[r, "d^1"]& A_{012}\\
& (a_0,a_1,a_2) \arrow[r, mapsto] & (a_0-a_1,a_2-a_0,a_1-a_2) &\\
&& (a_{01},a_{02},a_{12})\arrow[r, mapsto]& a_{01}+a_{02}+a_{12}\,.
\end{tikzcd}
\end{equation}
If $D=\sum_{y\subset X}n_y[y]$ is a divisor of $X$ we can define the subgroups
$$\mathbf A_X(D):=\prod_{y\subset X}\mathbb A^{(-n_y)}_y\,.$$
Note that $\mathbf A_X(D)$ is a well defined subgroup of $\mathbf A_X$ because $n_y=0$ for all but finitely many $y$. Let's define the subspaces 
$$A_{12}(D):=A_{012}\cap \mathbf A_X(D)=\mathbf A_X(D)\,.$$
$$A_1(D):=A_{01}\cap \mathbf A_X(D);\quad A_2(D):=A_{02}\cap \mathbf A_X(D);$$
in order to get the complex
\begin{equation}\label{ad_coomplexd}
\mathcal A_X(D):\quad A_0\oplus A_1(D)\oplus A_2(D)\xrightarrow{d^0_D} A_{01}\oplus A_{02}\oplus A_{12}(D)\xrightarrow{d^1_D} A_{012}
\end{equation}
such that the maps are the same of those in equation (\ref{ad_coomplex2}). Furthermore note that $\mathcal A_X=\mathcal A_X(0)$.
There is also the idelic version of complex (\ref{ad_coomplex2}):
\begin{equation}\label{id_coomplex}
\begin{tikzcd}
\mathcal A^\times_X: & A^\times_0\oplus A^\times_1\oplus A^\times_2 \arrow[r, "d_\times^0"]& A^\times_{01}\oplus A^\times_{02}\oplus A^\times_{12}\arrow[r, "d_\times^1"]& A^\times_{012}=\mathbf A_X^\times\\
& (a_0,a_1,a_2) \arrow[r, mapsto] & (a_0a^{-1}_1,a_2a^{-1}_0,a_1a^{-1}_2) &\\
&& (a_{01},a_{02},a_{12})\arrow[r, mapsto]& a_{01}a_{02}a_{12}\\
\end{tikzcd}
\end{equation}
and we have a well defined surjective map:

\begin{eqnarray*}
p:\;\; \ker(d^1_\times) &\to& \Div(X)\\
(\alpha,\beta,\alpha^{-1}\beta^{-1})&\mapsto & \sum_{y\subset X} v_y(\alpha_{x,y})[y]\,.
\end{eqnarray*}

\section{Idelic Deligne pairing}\label{sect2}
Let's still consider the arithmetic surface $\varphi: X\to B$, and let's denote with $\mathbf A_B$ and $\mathbf A_B^\times$ respectively the nonarchimedean parts of the one dimensional adeles and ideles associated to the base $B$ (in other words we are not considering the archimedean places of $K$). Remember that for any two divisors $D,E$ on $X$ with no common components we have:
$$i_x(D,E):= \len_{\mathscr O_{X,x}}\mathscr O_{X,x}/\left (\mathscr O_X(-D)_x+\mathscr O_X(-E)_x\right)\,.$$
The Deligne pairing is a bilinear and symmetric map:
$$
\left <\,,\,\right >\colon \Pic(X)\times\Pic(X)\to \Pic(B)
$$
which was introduced  for the first time in \cite{Deli}. More details about the construction of $\left <\,,\,\right >$ can be found in \cite[D.2.3]{dolce_phd}. First of all let's see how the Deligne pairing can be lifted to a pairing at the level of divisors on $X$ and with target in $\Pic(B)$:
\begin{proposition}\label{uni_pai_arith}
There exists a unique pairing
$$
[[\,,\,]]\colon \Div(X)\times\Div(X)\to \Pic(B)
$$
satisfying the following properties:
\begin{enumerate}
\item[$(1)$] It is bilinear and symmetric.
\item[$(2)$] It descends to the Deligne pairing 
$$
\left <\,,\,\right >\colon \Pic(X)\times\Pic(X)\to \Pic(B)\,.
$$
\item[$(3)$] If $D,E\in \Div(X)$ are two  divisors with no common components then $[[D,E]]$ is equal to the class  of the divisor $\left <D, E\right >$ in $\Pic(B)$ (here the bracket $\left <D,E\right >$ denotes the pushforward through $\varphi$ of the $0$-cycle collecting all the local intersection numbers  between $D$ and $E$). In other words:
$$[[D,E]]=\sum_{x\in D\cap E}[k(x): k(\varphi (x))]i_x(D,E)\,[\varphi(x)]\in\Pic(B)\,. $$ 
\end{enumerate}
\end{proposition}  
\proof
For any $D,E\in\Div(X)$ it is enough to put:
$$
[[D,E]]:=\left< \mathscr O_X(D), \mathscr O_X(E)\right>
$$
where on the right hand side we have the Deligne pairing between invertible sheaves. Uniqueness follows from properties $(1)$-$(3)$ and what is commonly called ``the moving lemma'' (\cite[Proposition 9.1.11]{Liu}).
\endproof
At this point we will try to work in complete analogy to the geometric case and we will use the Kato symbol defined in section \ref{subs1.1} to obtain the map, denoted below with a question mark, which makes the following diagram commutative:
\begin{equation}\label{id_arpai_dia}
\begin{tikzcd}
\ker(d^1_\times)\times \ker(d^1_\times)\arrow[rrdd, "{?}"]\arrow[d, two heads,"p\times p"] && \mathbf A^\times_B\arrow[dd ,two heads]\\
\Div(X)\times \Div(X)\arrow[rrd, "{[[\,,\,]]}",pos=.2]\arrow[d, two heads]  & &\\
\Pic(X)\times \Pic(X)\arrow[rr, "{\left <\,,\,\right >}"] &&\Pic(B)\cong\CH^1(B)\\
\end{tikzcd}
\end{equation}
As usual, fix a flag $x\in y$ with $\mathfrak z\in y(x)$ and assume that $\varphi(x)=b$, then we define
$$(\,,\,)_{x,\mathfrak z}:=(\,,\,)_{K_{x,\mathfrak z}|K_b}:K_{x,\mathfrak z}^\times\times K_{x,\mathfrak z}^\times\to K_b^\times
$$
where $(\,,\,)_{K_{x,\mathfrak z}|K_b}$ is the Kato symbol defined in section \ref{subs1.1}. Remember that depending on whether $y$ is horizontal or vertical, we have a different expression for $(\,,\,)_{x,\mathfrak z}$. Then we put:
$$
(\,,\,)_{x,y}:=\prod_{\mathfrak z\in y(x)}(\,,\,)_{x,\mathfrak z}:K^\times_{x,y}\times K^\times_{x,y}\to K^\times_b\,.
$$
It is important to point out that $(\,,\,)_{x,y}$ maps $\mathcal O^\times_{x,y}\times\mathcal O^\times_{x,y}$ to $\mathcal O^\times_b$.

\begin{proposition}\label{id_symb_ar}
The pairing $(\,,\,)_{x,y}$ is a skew-symmetric bilinear form on $K^\times_{x,y}$ satisfying the following properties:
\begin{enumerate}
\item[$(1)$]Let $r,s\in K^\times_x$, then for all but finitely many curves $y$ containing $x$ we have that $(r,s)_{x,y}=1$ and moreover $\prod_{y\ni x}(r,s)_{x,y}=1$.
\item[$(2)$] Let $y$ be a vertical curve and let $r,s\in K^\times_y$, then  $\prod_{x\in y}(r,s)_{x,y}=1$. In particular $(r,s)_{x,y}\in\mathcal O^\times_b$ for all but finitely many $x\in y$.
\end{enumerate}
\end{proposition}
\proof
Skew symmetry and bilinearity are clear. See \cite[Theorem 4.3]{otherLiu} for $(1)$; note that in \cite{otherLiu} the proof is made for $r,s\in K(X)^\times$, but it is easy to see that it actually works also for $r,s\in K^\times_x$. See \cite[Theorem 5.1]{otherLiu} for $(2)$.
\endproof

\begin{definition}\label{idelic_del}
The \emph{idelic Deligne pairing} 
$$\left <\,,\,\right >_i: \ker(d^1_\times)\times \ker(d^1_\times)\to \CH^1(B)$$ 
is given by:

\begin{equation}\label{idelic_del_eq}
(r,s)\mapsto \left <r,s\right >_i:=\sum_{b\in B} n_b(r, s)[b]\in\CH^1(B)
\end{equation}
such that:
\begin{equation}\label{idelic_del_prod}
n_b(r, s):= \sum_{\substack {x\in X_b,\\ y\ni x}}v_b\left((\gamma_{x,y},\beta_{x,y})_{x,y}\right)
\end{equation}
for $r=(\alpha,\beta, \alpha^{-1}\beta^{-1})$, $s=(\gamma,\delta,\gamma^{-1}\delta^{-1})\in \ker(d_1^{\times})$ and where $v_b$ is the complete discrete valuation on $K_b$. It is crucial to emphasize the fact the we consider $\sum_{b\in B} n_b(r, s)[b]$ in its linear equivalence class in $\CH^1(B)$ and not just as a divisor. By simplicity of notation we avoid to mention the canonical map $\Div(B)\to \CH^1(B)$.
\end{definition}
One can verify that definition \ref{idelic_del} makes sense:
\begin{proposition}
The summations (\ref{idelic_del_eq}) and (\ref{idelic_del_prod}) are finite.
\end{proposition}
\proof
Thanks to the second adelic restricted product over curves it is enough to check it only for a fixed nonsingular vertical curve $y\subset X_b$. Let's write $\beta_{x,y}=fs_{x,y}\in K(X)^\times\mathcal O^\times_x$, then 
\begin{equation}\label{finitessum}
(\gamma_{x,y},\beta_{x,y})_{x,y}=(\gamma_{x,y},f)_{x,y}(\gamma_{x,y},s_{x,y})_{x,y}\,.
\end{equation}
\cite[Lemma 5.2]{otherLiu} shows that $\prod_{x\in y}(\gamma_{x,y},f)_{x,y}$ is convergent, and this means that $(\gamma_{x,y},f)_{x,y}\in\mathcal O^\times_b$ for all but finitely many $x\in y$. Moreover if $p=\cha k(b)$, then $p$ is a uniformizing parameter for $K_{x,y}$, so $\gamma_{x,y}=p^r c_{x,y}$ with $c_{x,y}\in\mathcal O^\times_{x,y}$. Thus:

$$(\gamma_{x,y},s_{x,y})_{x,y}=(p^r,s_{x,y})_{x,y}(c_{x,y},s_{x,y})_{x,y}$$
Obviously $(c_{x,y},s_{x,y})_{x,y}\in\mathcal O^\times_b$, so in order to finish the proof it remains to show that $(p^r,s_{x,y})_{x,y}$ lies in $\mathcal O^\times_b$ too. Just for simplicity of calculations  let's assume that $K_{x,y}$ is a standard field and  $K_{x,y}=K_p\{\!\{t\}\!\}$ (the argument works easily also for non standard fields). By the explicit expression of the Kato's residue homomorphism (cf. \cite[equation (8)]{otherLiu}) we can calculate that:

$$(p^r, s_{x,y})_{x,y}= N_{K_p|K_b}\left(p^{rw}\right)$$
where $w\in\mathbb Z$ is the  winding number associated to $s_{x,y}$ (see \cite[equation (7)]{otherLiu} for details). But we know that $s_{x,y}\in\mathcal O^\times_x=\mathcal O_{K_p}^\times+t\mathcal O_{K_p}[[t]]$, thus:
$$s_{x,y}=a+t\sum_{i\ge 0}a_it^i=a\left(1+t\sum_{i\ge 0}\frac{a_i}{a}t^i\right)\,.$$
It follows that $w=0$ and the proof is complete. 
\endproof

\begin{remark}\label{idelic_dec}
For any $b\in B$ we have the following decomposition for the big product  (\ref{idelic_del_prod}):
$$
\sum_{\substack {x\in X_b,\\ y\ni x}}v_b\left((\,,\,)_{x,y}\right)=\sum_{\substack {y\subset X_b,\\ x\in y}} v_b\left((\,,\,)_{x,y}\right)   +\sum_{\substack {x\in X_b,\\ y\ni x,\\ y\text{ horiz.}}}v_b\left((\,,\,)_{x,y}\right)\,.
$$
\end{remark}
We put $\left <\,,\, \right>_i$ as the undetermined function in diagram (\ref{id_arpai_dia}) and we have the following fundamental result:
\begin{theorem}\label{idelic_big_teo}
Consider the notation of diagram (\ref{id_arpai_dia}). The pairing $\left<\;,\;\right>_i$ satisfies the following properties:
\begin{enumerate}
\item[$(1)$] It is bilinear and symmetric.
\item[$(2)$] Let $r,s,r',s'\in\ker(d^1_{\times})$ such that $p(r)=p(r')$ and $p(s)=p(s')$, then $\left<r,s\right>_i=\left<r',s'\right>_i$.
\item[$(3)$] It descends naturally to a pairing $H^1(\mathcal A_X^\times)\times H^1(\mathcal A_X^\times)\to \Pic(B)$.
\end{enumerate}
\end{theorem}
\proof
Let's fix $r=(\alpha,\beta,\alpha^{-1}\beta^{-1})$, $s=(\gamma,\delta,\gamma^{-1}\delta^{-1})\in\ker(d^1_{\times})$; moreover we can fix $b\in B$ and work componentwise.
\\
$(1)$ Bilinearity is clear. We will show that as elements of $\Div(B)$ we have $\left<r,s\right>_i=\left<s,r\right>_i+(f)$ with $f\in K^\times$. 
 For any flag $x\in y$: $\alpha_{x,y}^{-1}\beta_{x,y}^{-1}\in\mathcal O^\times_{x,y}$ and $\gamma_{x,y}^{-1}\delta_{x,y}^{-1}\in\mathcal O^\times_{x,y}$ so  we have that:
\begin{equation}\label{difficult_eq}
\begin{gathered}
0=\sum_{\substack {x\in X_b,\\ y\ni x}}v_b\left((\alpha_{x,y}^{-1}\beta_{x,y}^{-1},\gamma_{x,y}^{-1}\delta_{x,y}^{-1})_{x,y}\right)=\\
=\underbrace{\sum_{\substack {x\in X_b,\\ y\ni x}}v_b\left((\alpha_{x,y},\gamma_{x,y})_{x,y}\right)}_{(i)}+\underbrace{\sum_{\substack {x\in X_b,\\ y\ni x}}v_b((\alpha_{x,y},\delta_{x,y})_{x,y})}_{(ii)}+\\
+\underbrace{\sum_{\substack {x\in X_b,\\ y\ni x}}v_b((\beta_{x,y},\gamma_{x,y})_{x,y})}_{(iii)}+\underbrace{\sum_{\substack {x\in X_b,\\ y\ni x}}v_b((\beta_{x,y},\delta_{x,y})_{x,y})}_{(iv)}\,.\end{gathered}
\end{equation}
Now we analyze in detail the underbraced terms in equation (\ref{difficult_eq}): for $(i)$ we have the following decomposition thanks to remark \ref{idelic_dec}:
$$(i)=\sum_{\substack {y\subset X_b,\\ x\in y}}v_b((\alpha_{x,y},\gamma_{x,y})_{x,y})  + \sum_{\substack {x\in X_b,\\ y\ni x,\\ y\text{ horiz.}}}v_b((\alpha_{x,y},\gamma_{x,y})_{x,y})=$$
$$
^{(\text{prop.} \ref{id_symb_ar}(2))}=0+\sum_{x\in X_b} \sum_{\substack {y\ni x,\\ y\text{ horiz.}}}v_b((\alpha_{x,y},\gamma_{x,y})_{x,y})\,.
$$
By definition we have that $(ii)=n_b(s, r)$ and $(iii)=-n_b(r, s)$. Finally:
$$(iv)=\sum_{x\in X_b}\sum_{y\ni x}(\beta_{x,y},\delta_{x,y})_{x,y}=^{(\text{prop.} \ref{id_symb_ar}(1))}\;0$$
By substituting in equation (\ref{difficult_eq}) we conclude that  
\begin{equation}\label{lambdas}
n_b(r, s)=n_b(s, r)+\sum_{\substack {x\in X_b,\\ y\ni x,\\ y\text{ horiz.}}}v_b((\alpha_{x,y},\gamma_{x,y})_{x,y})\,.
\end{equation}
Let $y$ be an horizontal curve and let $x\in y$ such that $\varphi(x)=b$, then the coefficient field of $K_{x,\mathfrak z}$ is $k(y)_x$. The two dimensional valuation $v_{x,\mathfrak z}$ extends the valuation $v_y$ on $k(y)$ and moreover that the norm $N_{k(y)_x|K_b}$ extends $N_{k(y)|K}$. It follows that  $(\,,\,)_{x,\mathfrak z}$ extends the one dimensional tame symbol 
$$(\,,\,)_{y}:=N_{k(y)|K}\circ (\,,\,)_{k(y)}:K^\times_y\times K^\times_y\to k(y)^\times\to K^\times\,.$$
This means that for any two elements $u,v\in K_y$, where $y$ is horizontal, we have that:
$$(u,v)_{x,y}=(u,v)_y\in K$$
for any $x\in y$. Therefore we can rewrite equation (\ref{lambdas}):
\begin{equation}\label{lambdass}
n_b(r, s)=n_b(s, r)+\sum_{y\text{ horiz.}}v_b((\alpha_{x,y},\gamma_{x,y})_{y})\,.
\end{equation}
Let's put $f=\prod_{y\text{ horiz.}}(\alpha_{x,y},\gamma_{x,y})_{y}\in K^\times$, then equation (\ref{lambdass}) implies the following equality:
$$\left<r,s\right>_i=\left<s,r\right>_i+\sum_{b\in B} v_b(f)[b]=\left<s,r\right>_i+(f)\,.$$

$(2)$ Let  $r'=(\alpha',\beta', (\alpha')^{-1}(\beta')^{-1})$  and $s'=(\gamma',\delta', (\gamma')^{-1}(\delta')^{-1})$. Since $p(r)=p(r')$ and $p(s)=p(s')$, then $v_y(\alpha_{x,y})=v_y(\alpha'_{x,y})$ and $v_y(\gamma_{x,y})=v_y(\gamma'_{x,y})$. This means that $\gamma'_{x,y}=f_{x,y}\gamma_{x,y}$ and $\alpha'_{x,y}=g_{x,y}\alpha_{x,y}$ for $f_{x,y},g_{x,y}\in\mathcal O^\times_y$ (for any $x\in y$). Then we have the following chain of equalities depending on what we showed in claim $(1)$:

$$n_b(r',s')=\sum_{\substack {x\in X_b,\\ y\ni x}}v_b((f_{x,y}\gamma_{x,y},\beta'_{x,y})_{x,y})=$$

$$=\sum_{\substack {x\in X_b,\\ y\ni x}}((f_{x,y},\beta'_{x,y})_{x,y})+\sum_{\substack {x\in X_b,\\ y\ni x}}v_b((\gamma_{x,y},\beta'_{x,y})_{x,y})=$$

$$=\sum_{\substack {x\in X_b,\\ y\ni x}}v_b((f_{x,y},\beta'_{x,y})_{x,y})+n_b(r',s)=\sum_{\substack {x\in X_b,\\ y\ni x}}v_b((f_{x,y},\beta'_{x,y})_{x,y})+n_b(s,r')+v_b(f)=$$
$$=\sum_{\substack {x\in X_b,\\ y\ni x}}v_b((f_{x,y},\beta'_{x,y})_{x,y})+\sum_{\substack {x\in X_b,\\ y\ni x}}v_b((g_{x,y}\alpha_{x,y},\delta_{x,y})_{x,y})+v_b(f)=(\ast)$$
Where $f\in K^\times$.
$$
(\ast)=\sum_{\substack {x\in X_b,\\ y\ni x}}v_b((f_{x,y},\beta'_{x,y})_{x,y})+ \sum_{\substack {y\subset X_b,\\ x\in y}}v_b((g_{x,y},\delta_{x,y})_{x,y})+\sum_{\substack {x\in X_b,\\ y\ni x}}v_b((\alpha_{x,y},\delta_{x,y})_{x,y})+v_b(f)=$$
$$=\underbrace{\sum_{\substack {x\in X_b,\\ y\ni x}}v_b((f_{x,y},\beta'_{x,y})_{x,y})}_{(i)}+ \underbrace{\sum_{\substack {x\in X_b,\\ y\ni x}}v_b((g_{x,y},\delta_{x,y})_{x,y})}_{(ii)}+n_b(r,s)+v_b(fg)\,.
$$
Note that in the last line we used the fact that $n_b(s,r)=n_b(r,s)+v_b(g)$ for $g\in K^\times$. We have to show that the terms $(i)$ and $(ii)$ are valuations at $b$ of elements of $K^\times$. Since $(\alpha')^{-1}_{x,y}(\beta')^{-1}_{x,y}\in\mathcal O^\times_{x,y}$, we have:
$$0=\sum_{\substack {x\in X_b,\\ y\ni x}}v_b((f_{x,y},(\alpha')^{-1}_{x,y}(\beta')^{-1}_{x,y})_{x,y})=
$$
$$
=\sum_{\substack {x\in X_b,\\ y\ni x}}v_b((f_{x,y},\alpha'_{x,y})_{x,y})+\sum_{\substack {x\in X_b,\\ y\ni x}}v_b((f_{x,y},\beta'_{x,y})_{x,y})=
$$
$$
=v_b(h)+\sum_{\substack {y\subset X_b,\\ x\in y}}v_b((f_{x,y},\beta'_{x,y})_{x,y})\,.$$
with $h=\prod_{y\text{ horiz.}}(f_{x,y},\alpha'_{x,y})_{y}\in K^\times$. For $(ii)$ the argument is similar.\\

$(3)$ Let $r,s\in\im(d^0_\times)$. It means that $\alpha=lm^{-1}$, $\beta=tl^{-1}$, $\gamma=uv^{-1}$, $\delta=zu^{-1}$ for $l,u\in A^\times_0=K(X)^\times$,  $m,v\in A^\times_1$ and $t,z\in A_2^\times$. So:

\begin{equation}\label{difficult_eq1}
\begin{gathered}
n_b(r,s)=\sum_{\substack {x\in X_b,\\ y\ni x}}v_b((uv^{-1}_{x,y},t_{x,y}l^{-1})_{x,y})=\\
=\sum_{\substack {x\in X_b,\\ y\ni x}}v_b((u,t_{x,y})_{x,y})+\sum_{\substack {x\in X_b,\\ y\ni x}}v_b((u,l^{-1})_{x,y})+\\
+\sum_{\substack {x\in X_b,\\ y\ni x}}v_b(v^{-1}_{x,y},t_{x,y})_{x,y})+\sum_{\substack {x\in X_b,\\ y\ni x}}v_b((v^{-1}_{x,y}, l^{-1})_{x,y})\,.
\end{gathered}
\end{equation}
Now it is enough to appeal to one of the arguments previously used to conclude that each summand of equation (\ref{difficult_eq1}) is either $0$ or of the form $v_b(f)$ for $f\in K^\times$. It means that $\left<r,s\right>_i=0$ in $\CH^1(B)$.
\endproof
We want to give an alternative formula for the coefficient $n_b(r,s)$. Notice that:
\begin{equation}
\begin{gathered}
-\sum_{\substack {x\in X_b,\\ y\ni x}}v_b((\alpha_{x,y},\gamma^{-1}_{x,y}\delta^{-1}_{x,y})_{x,y})=\\
=\sum_{\substack {x\in X_b,\\ y\ni x}}v_b((\alpha_{x,y},\gamma_{x,y})_{x,y})+\sum_{\substack {x\in X_b,\\ y\ni x}}v_b((\alpha_{x,y},\delta_{x,y})_{x,y})=\\
=v_b(f)+\sum_{\substack {x\in X_b,\\ y\ni x}}v_b((\alpha_{x,y},\delta_{x,y})_{x,y})=v_b(f)+n_b(s,r)=v_b(fg)+n_b(r,s)
\end{gathered}
\end{equation}
for $f,g\in K^\times$. Therefore,  we can also express:

\begin{equation}
n_b(r,s)=-\sum_{\substack {x\in X_b,\\ y\ni x}}v_b((\alpha_{x,y},\gamma^{-1}_{x,y}\delta^{-1}_{x,y})_{x,y})\,.
\end{equation}
In particular if $y$ is a horizontal curve:
\begin{equation}\label{alt_pairingg1}
-v_b((\alpha_{x,y},\gamma^{-1}_{x,y}\delta^{-1}_{x,y})_{x,y})=v_y(\alpha_{x,y})v_b\left(N_{k(y)_x|K_b}\left(\overline{\gamma^{-1}_{x,y}\delta^{-1}_{x,y}}\right)\right)\,.
\end{equation}
The following lemmas are fundamental in order to understand the relationship between $\left<r,s\right>_i$ and Deligne pairing.
\begin{lemma}\label{algeolemma}
Let $X_b\subset X$ the fiber over $b\in B$ and assume that $X_b$ has at least two irreducible components. If $D\subset X_b$ is an integral curve, then there exists a divisor $D'\sim D$ such that $D'$ doesn't have components contained in $X_b$.
\end{lemma}
\proof
Consider $\Gamma$ running amongst all irreducible components of $X_b$, then put
$$S:=\bigcup_{\substack{\Gamma\subset X_b\\ \Gamma\neq D}}(\Gamma\cap D)\,.$$
By the moving lemma we can find $D'\sim D$ not passing by $S$. It is clear by the definition of $S$ that $D'$ cannot have vertical components contained in $X_b$.
\endproof

\begin{lemma}\label{int_val}
Let $D,E$ two prime divisors on $X$ and let $x\in D\cap E$ a nonsingular point for both $D$ and $E$. Moreover let $d_x,e_x\in\mathscr O_{X,x}$ be the local equations at $x$ of $D$ and $E$ respectively. Then we have the equality:
$$v^{(1)}_{x,D}(\overline{e_x})=i_x(D,E)$$
where $ v^{(1)}_{x,D}:E^\times_{x,D}\to\mathbb Z$ is the one dimensional valuation and $\overline{e_x}\in E^\times_{x,D}$ is the natural projection through the map $\mathcal O_{x,D}\to E^\times_{x,D}$.
\end{lemma} 
\proof
Put  $y=D$ and $v=v^{(1)}_{x,D}$. First of all notice that $\mathscr O_{X,x}\subseteq\mathscr O_{X,y}$, therefore $\overline{e_x}\in \mathscr O_{y,x}$  and it is the image on the natural map $\mathscr O_{X,x}\to\mathscr O_{y,x}\subset k(y)$. We have to show that $v(\overline{e_x})=\len_{\mathscr O_{X,x}}\frac{\mathscr O_{X,x}}{(d_x,e_x)}$, but we know that $\mathscr O_{y,x}=\frac{\mathscr O_{X,x}}{d_x\mathscr O_{X,x}}$, thus
$$\len_{\mathscr O_{X,x}}\frac{\mathscr O_{X,x}}{(d_x,e_x)}=\len_{\mathscr O_{y,x}}\frac{\mathscr O_{y,x}}{\overline{e_x}\mathscr O_{y,x}}=v(\overline{e_x})\,.$$
\endproof
\begin{theorem}\label{arith_id_pai1}
If $r,s\in \ker(d^1_\times)$ such that $D=p(r)$ and $E=p(s)$ are two nonsingular prime divisors on $X$ with no common components, then $\left<r,s\right>_i=[[D,E]]$.
\end{theorem}
\proof 
Fix $r=(\alpha,\beta, \alpha^{-1}\beta^{-1})$, $s=(\gamma,\delta,\gamma^{-1}\delta^{-1})$. We want to show that it is enough to restrict to the case when either $D$ or $E$ is horizontal. In any case, by theorem \ref{idelic_big_teo}(2) we always choose $\delta_{x,y}$ in the following way: $\delta_{x,y}=1$ if $x\notin D\cap E$ and $\delta_{x,y}=t_x^{-1}$, where $t_x\in\mathscr O_{X,x}$ is the local equation of $E$ at $x$, if $x\in D\cap E$. For any $y\neq D$, $\alpha_{x,y}\in\mathcal O^\times_y$, since $p(a)=D$, therefore:
\begin{equation}\label{eqnb}
n_b(r,s)=-\sum_{x\in D\cap X_b}v_b((\alpha_{x,D},\gamma^{-1}_{x,D}\delta^{-1}_{x,D})_{x,D})\,.
\end{equation}
If $D\subseteq X_b$ and $E\subseteq X_{b'}$ with $b\neq b'$, then by proposition \ref{id_symb_ar}(2) and the choice of $\delta_{x,y}$ we have:
 \begin{equation}
n_b(r,s)=-\sum_{x\in D\cap X_b}v_b((\alpha_{x,D},\gamma^{-1}_{x,D})_{x,D})=0\,.
\end{equation}
So in such a particular case $\left<r,s\right>_i=[[D,E]]=0$.\\
If $D,E\in X_b$ we can apply lemma \ref{algeolemma} and find a divisor $D'=\sum_j n_j\Gamma_j\sim D$ such that $\Gamma_j\not\subset X_b$. Clearly

$$[[\Gamma_j,E]]=\sum_j n_j[[\Gamma_j,E]]$$ 
therefore from now on we can restrict our calculation to the case where either $D$ or $E$ is horizontal. By symmetry we can fix $D$ to be horizontal and we denote with $K(D)$ its function field. In this case we have an explicit expression given by equation (\ref{alt_pairingg1}):
$$
n_b(r,s)=\sum_{x\in D\cap X_b}v_b\left(N_{K(D)_x|K_b}\left(\overline{\gamma^{-1}_{x,D}t_x}\right)\right)=\sum_{x\in D\cap X_b}v_b\left(N_{K(D)|K}\left(\overline{\gamma^{-1}_{x,D}t_x}\right)\right)=
$$
$$
=\sum_{x\in D\cap X_b}v_b\left(N_{K(D)_x|K_b}\left(\overline{\gamma^{-1}_{x,D}}\right)\right)+\sum_{x\in D\cap E\cap X_b}v_b\left(N_{K(D)_x|K_b}\left(\overline{t_x}\right)\right)\,.
$$
Now by the theory of extensions of valuation fields (see \cite[II(2.5)]{fesbook}), we know that if $v_x:=v^{(1)}_{x,D}$ is the valuation on $K(D)_x$, then:
$$v_x=\frac{1}{[k(x):k(b)]}v_b\circ N_{K(D)_x|K_b}\,.$$
Therefore we obtain:
\begin{equation}\label{last_fund_rel}
n_b(r,s)=\sum_{x\in D\cap X_b}[k(x):k(b)]v_x\left(\overline{\gamma^{-1}_{x,D}}\right)+\sum_{x\in D\cap E\cap X_b}[k(x):k(b)]v_x\left(\overline{t_x}\right)\,.
\end{equation}
Put by simplicity $f=\overline{\gamma^{-1}_{x,D}}\in K(D)^\times$, consider the restricted morphism of arithmetic curves $\varphi:D\to B$ and the principal divisor $(f)\in \Princ(D)$, then:
$$\varphi_\ast((f))=\sum_{b\in B}\left(\sum_{x\in D\cap X_b} [k(x):k(b)]v_x(f)\right)[b]\,.$$
Moreover $v_x\left(\overline{t_x}\right)=i_x(D,E)$ by lemma \ref{int_val}. Equation \ref{last_fund_rel} implies that in $\Div(B)$ we have the following equality:
$$\left<r,s\right>_i=\varphi_\ast((f))+[[D,E]]\,.$$
But by \cite[7 Remark 2.19]{Liu} we know that $\varphi_\ast((f))=\left(N_{K(D)|K}(f)\right)\in\Princ(B)$, so the proof is complete. 
\endproof
We obtained the idelic representation of Deligne pairing:
\begin{corollary}
Diagram (\ref{id_arpai_dia}) is commutative.
\end{corollary}
\proof
For any two divisors $D,E\in \Div(X)$ define the pairing:
$$\Theta(D,E):=\left<r',s'\right >_i$$
for a choice  of $r',s'\in\ker(d^1_\times)$ such that $p(r')=D$ and $p(s')=E$. By theorem \ref{idelic_big_teo}(2) $\Theta$ is well defined and moreover by \ref{idelic_big_teo}(1), \ref{idelic_big_teo}(3) and \ref{arith_id_pai1} we can conclude that $\Theta(D,E)=D.E$. Thus for any $a,b\in\ker(d^1_\times)$ we have that:
$$\left <r,s\right>_i=\Theta(p(r),p(s))=p(r).p(s)\,.$$
\endproof

\section{Adelic Deligne pairing}\label{sect3}
A clever and quick adelic interpretation of  intersection theory on algebraic surfaces is given in \cite{fe0} and the strategy is very simple: first of all one defines the adelic Euler-Poincare characteristic $\chi_{a}(\cdot)$ which associates an integer to any divisor on the surface (or more in general to any invertible sheaf) by using just data coming from the adelic complex. Then the adelic intersection pairing is defined accordingly to equation (\ref{geom_pai}) by using $\chi_{a}(\cdot)$ instead of the usual Euler-Poincare characteristic.  Here we try to follow the same approach; so, it is evident that we have to  define the adelic determinant of the cohomology
 $$\det{\!}_a R\varphi_\ast:\Pic(X)\to\Pic(B)$$ 
which should be a function involving only adelic data, and then the adelic Deligne pairing according to theorem  \ref{det_coh_deli}.

For any coherent sheaf $\mathscr F$ on $B$ and any closed point $b\in B$, we define the following  objects:
$$K_b(\mathscr F):=\mathscr F_\xi\otimes_K K_b=(\mathscr F_b\otimes_{\mathscr O_{B,b}} K)\otimes_K K_b=\mathscr F_b\otimes_{\mathscr O_{B,b}} K_b\,,$$
$$\mathcal O_b(\mathscr F):=\mathscr F_b\otimes_{\mathscr O_{B,b}} \mathcal O_b\,.$$
$$\mathbf A_B(\mathscr F):=\sideset{}{'}\prod_{b\in B} K_b(\mathscr F)$$
where the restricted product is taken with respect to the rings $\mathcal O_b(\mathscr F)$, and
$$\mathbf A_B(\mathscr F)(0):=\prod_{b\in B} \mathcal O_b(\mathscr F)\,.$$
Moreover recall that we have the following one dimensional adelic complex given by:

\begin{eqnarray*}\label{adelic_complexx}
\mathcal A_B(\mathscr F):\quad 0\rightarrow \mathscr F_{\xi}\oplus\mathbf A_B(\mathscr F)(0)&\rightarrow &\mathbf A_B(\mathscr F)\rightarrow 0 \\
(f,(\alpha_b)_b) &\mapsto& (f-\alpha_b)_b
\end{eqnarray*}
It is important to point out that we want to  consider $\mathcal A_B(\mathscr F)$ as a complex of $O_K$-modules in the natural way.
\begin{definition}
Let $D$ be a divisor on $X$ satisfying proposition \ref{fund_det_coho}. For any invertible sheaf $\mathscr L$ on $B$ we put by simplicity $\mathscr G:=\varphi_\ast\mathscr L(D)$ and $\mathscr H:=\varphi_\ast(\mathscr L(D)/\mathscr L)$. Then the \emph{adelic determinant of cohomology} is given by:

$$\det{\!}_{a} R\varphi_\ast(\mathscr L):=\det \left(H^0(\mathcal A_B(\mathscr G))\right)\otimes \left(\det\left(H^0(\mathcal A_B(\mathscr H))\right)\right)^{\ast}$$
Where with $\ast$ we denote the algebraic dual. Note that $\det{\!}_{a} R\varphi_\ast(\mathscr L)$ is a $O_K$-module, but by abuse of notation we can consider it as an element in $\Pic(B)$ after taking the associated sheaf $(\det{\!}_{a} R\varphi_\ast(\mathscr L))^\sim$ (see \cite[5.1.2]{Liu}). In other words we omit the operator $\sim$ by simplicity of notations. 

\end{definition}
\begin{definition}
The \emph{adelic Deligne pairing}  between two invertible sheaves $\mathscr L$ and $\mathscr M$ on $B$ is defined as:
$$\left<\mathscr L,\mathscr M\right>_ a:=$$
$$
\det{\!}_{a}R\varphi_\ast(\mathscr O_X)\otimes (\det{\!}_{a} R\varphi_\ast(\mathscr L))^{-1}\otimes(\det{\!}_{a} R\varphi_\ast(\mathscr M))^{-1}\otimes\det{\!}_{a} R\varphi_\ast(\mathscr L\otimes \mathscr M)
$$
\end{definition}
It is immediate to verify that the definition of the adelic Deligne pairing coincides with the usual Deligne pairing by using equation (\ref{det_coh_deli}). Indeed thanks to \cite{huber} $H^0(\mathcal A_B(\mathscr G))\cong H^0(B,\mathscr G)$ and $H^0(B,\mathscr G)^\sim\cong\mathscr G$ by affine Serre's theorem (see \cite[II,4]{Fac}). Obviously the same holds for $\mathscr F$.

\begin{appendices}
\section{Topics in $K$-theory}\label{k_th}
Algebraic $K$-theory is a very wide subject with a long history. It can be approached in many different ways and several links can be build between all approaches (see for example \cite{srini}). This appendix  is not a short introduction to algebraic $K$-theory, but just a mere collection of definition and notations needed in this text.

\begin{definition}
Let $G$ be an abelian group, and fix an integer $\ge 1$. A \emph{$r$-Steinberg map} is an homomorphism of $\mathbb Z$-modules $f: (F^\times)^{\oplus_r}\to G$ such that $f(a_1,\ldots,a_r)=0$ whenever there exist two indexes $i,j$ such that $i\neq j$ and $a_i+a_j=1$.
\end{definition}
Let's denote with $\catname{St}(r)$ the category whose objects are the $r$-Steinberg maps  $f: (F^\times)^{\oplus_r}\to G$ and the morphisms are the commutative diagrams:

$$
\begin{tikzcd}
(F^\times)^{\oplus_r}\arrow[r,"f"]\arrow[d, "g"']& G \\
H\arrow[ur, "\phi"']&  
\end{tikzcd}
$$
where $\phi$ is a group homomorphism.
\begin{proposition}\label{K-gr}
The category $\catname{St}(r)$ has the initial object.
\end{proposition}
\proof
We construct the initial objects by hands. Let's define
$$K_r(F):=\underbrace{F^\times\otimes_{\mathbb Z}\ldots\otimes _{\mathbb Z}F^\times}_{\text{$r$ times}}\,\big /\, S$$
where $S$ is the (multiplicative) subgroup generated by the following set:
$$\{a_1\otimes\ldots\otimes a_r\colon a_i+a_j=1 \text{ for some } i\neq j\}\,.$$
The natural image of a pure tensor $a_1\otimes\ldots\otimes a_r$ in $K_r(F)$ is denoted by $\left\{a_1,\ldots,a_r\right\}$.  Clearly we have an induced map:
\begin{eqnarray*}
\{\;\}:\;\;(F^\times)^{\oplus_r}&\to&K_r(F)\\
(a_1,\ldots,a_r)&\mapsto& \left\{a_1,\ldots,a_r\right\}
\end{eqnarray*}
At this point it is straightforward to see that $\{\;\}:(F^\times)^{\oplus_r}\to K_r(F)$ is the initial object for $\catname{St}(r)$.
\endproof
\begin{definition}\label{miln_k}
For $r=0$ we put $K_0(F):=\mathbb Z$ and in general we call the group $K_r(F)$ constructed in proposition \ref{K-gr} \emph{the $r$-th $K$-group of $F$}. Note that $K_1(F)=F^\times$. The map  $\{\;\}:(F^\times)^{\oplus_r}\to K_r(F)$ is called the $r$-th symbol map and in the cases $r=0,1$ it is just the identity.
\end{definition}

\begin{remark}
The groups introduced in definitons \ref{miln_k} are usually called Milnor $K$-groups and the standard notation is $K_r^{M}$. However in this text we can simplify the notation.
\end{remark}
The construction $K_r(\;)$  is functorial,  in fact let $f:F^\times\to L^\times$ be a group homomorphism, then the composition:
$$
(F^\times)^{\oplus r}\xrightarrow{f^{\oplus r}} (L^\times)^{\oplus r}\xrightarrow{\{\,,\,\}} K_r(L)
$$
is evidently a Steinberg map. By the universal property it induces a morphism $K_r(f): K_r(F)\to K_r(L)$.

When $F$ is a complete discrete valuation field there exists a nice relationship between $K$-groups of $F$ and $K$-groups of the residue fields:
\begin{theorem}\label{Miln_boundary}
Let $F$ be a discrete valuation field (not necessarily complete) then there is a unique group homomorphism:
$$\partial_r: K_r(F)\to K_{r-1}(\overline F)$$
satisfying the following property:
$$\partial_r(\{x_1,\ldots,x_{r-1},\varpi\})=\{\overline{x_1},\ldots,\overline{x_{r-1}}\}$$
for any local parameter $\varpi$ of $F$ and any $x_1,\ldots,x_{r-1}\in \mathcal O^\times_F$\,.
\end{theorem}
\proof See \cite{miln}.
\endproof

\begin{definition}
The map $\partial_r$ described in theorem \ref{Miln_boundary} is called \emph{the (Milnor) $r$-th boundary map}.
\end{definition}

Consider the \emph{tame symbol} for a complete discrete valuation field $(F,v)$:
\begin{equation}\label{tame_sym}
\begin{aligned}
(\,,\,)_F:F^\times\times F^\times &\to {\overline F}^\times\\
(a,b)&\mapsto (a,b)_F=(-1)^{v(a)v(b)}\,\overline{a^{v(b)}b^{-v(a)}}\,. 
\end{aligned}
\end{equation}
We have a nice description of the boundary map $\partial_2$ in relation to the tame symbol.  By the universal property of $K_2(F)$, the tame symbol $(\,,\,)_F$ induces a unique map $\Psi:K_2(F)\to\overline F^\times=K_1(F)$ such that $\Psi(\{\,,\,\})=(\,,\,)_F$. Let $a\in \mathcal O^\times_F$ and let $\varpi$ be  a local parameter for $F$, then $\Psi(\{a,\varpi\})=(a,\varpi)_F=\overline a$; this actually means that $\partial_2=\Psi$. In other words the $2$-nd boundary map for a complete discrete valuation field is exactly the map induced naturally by the tame symbol.

For a discrete valuation field $F$ (not necessarily complete) we have the multiplicative group $U_F^{(i)}:=1+\mathfrak p_F^i$  for $i\ge 1$ and we have also the $K$-theoretic version of it:
$$U^iK_r(F):=\{\{a_1\ldots a_r\}\in K_r(F)\colon a_j\in U_F^{(i)} \; \forall j=1,\ldots,r\}$$
and we put:
\begin{equation}\label{completedK}
\widehat{K}_r(F):=\varprojlim_i K_r(F)/U^iK_r(F)\,.
\end{equation}
Clearly we have a natural homomorphism $K_r(F)\to\widehat{K}_r(F)$ and moreover if $\widehat F$ is the completion of $F$ there is an isomorphism $\widehat{K}_r(F)\cong\widehat{K}_r(\widehat F)$. Now put $L=\fr\left(\mathcal O_F[[t]]\right )$, for any prime ideal $\mathfrak p$ of height $1$ in $\mathcal O_F[[t]]$ we have that  $\mathcal O_F[[t]]_{\mathfrak p}$ is a discrete valuation ring and in particular $F\{\!\{t\}\!\}$ is the completion of $L$ at $\mathfrak p= \mathfrak p_F\mathcal O_F[[t]]$. Consider the set: 
$$\mathfrak S:=\{\mathfrak p\in\spec(\mathcal O_F[[t]])\colon \hei \mathfrak p=1, \mathfrak p \neq \mathfrak p_F\mathcal O_F[[t]]\}\,,$$
and for any $\mathfrak p\in \mathfrak S$ let's denote with $\partial^{(\mathfrak p)}_r:K_r(L)\to K_{r-1}(k(\mathfrak p))$ the $r$-th boundary map relative to the valuation defined by $\mathfrak p$. 
\begin{definition}
For $r\ge 1$, the $r$-th (Kato) residue map on $L$ is given by the following composition:
$$
\begin{tikzcd}
\displaystyle\res_{L}^{(r)}: K_r(L)\arrow[rrr,"(\partial^{(\mathfrak p)}_r)_{\mathfrak p\in\mathfrak S}"] & & &\displaystyle \bigoplus_{\mathfrak p\in\mathfrak S} K_{r-1} (k(\mathfrak p))\arrow[rrr, " \sum_{\mathfrak p\in\mathfrak S} K_r(N_{k(\mathfrak p)|F})"] & & &  K_{r-1}(F)
\end{tikzcd}
$$
\end{definition}
\begin{theorem}\label{Kato_res}
The $r$-th residue map satisfies:
$$\res_{L}^{(r)}\left(U^{(i)}K_r(L)\right)\subseteq U^{(i)} K_{r-1}(F)\quad \forall i\ge 1 $$
therefore it induces a homomorphism:
$$\res_{F\{\!\{t\}\!\}}^{(r)}:\widehat {K}_r(F\{\!\{t\}\!\})\cong \widehat {K}_r(L)\to \widehat{K}_{r-1}(F)\,.$$
\end{theorem}
\proof
See \cite[Theorem 1]{kato}.
\endproof

\section{Determinant of cohomology}\label{det_coho}
For an algebraic surface $Z$ over a field $k$, intersection theory can be introduced by using the Euler-Poincare characteristic $\chi_k:\Coh(Z) \to \mathbb Z$ ``restricted'' to $\Pic(Z)$. In fact, the intersection number between two invertible sheaves $\mathscr L$ and $\mathscr M$ on $Z$ can be calculated by the following formula:
\begin{equation}\label{geom_pai}
\mathscr L.\,\mathscr M:= \chi_k(\mathscr O_Z)-\chi_k(\mathscr L^{-1})-\chi_k(\mathscr M^{-1})+\chi_k(\mathscr L^{-1}\otimes \mathscr M^{-1})
\end{equation}
 In the Arakelov setting given by the arithmetic surface $\varphi:X\to B$, we want to define a map $\Coh(X) \to \Pic(B)$ such that, when we take the ``restriction'' to $\Pic(X)$, we  obtain a formula, similar to (\ref{geom_pai}), relating our map to the Deligne pairing. In other words, we would like to have the arithmetic equivalent notion of the Euler-Poincare characteristic. The answer to our query will be the determinant of the cohomology, denoted by $\det R\varphi_\ast$, and in this section we are going to construct it step by step.
 
The first thing to notice is that $\chi_k$ is a cohomological object and in the case of $\varphi: X\to B$ the ``relative cohomology" is captured by the higher direct image functors $R^i\varphi_\ast$. By keeping in mind that the output of the determinant of the cohomology should be an invertible sheaf on the base $B$, in analogy with the definition of $\chi_k$, the most reasonable definition should be something like:
\begin{equation}\label{wrongdet_coho}
\det R\varphi_\ast(\mathscr F):=\bigotimes_{j\ge 0} (\det R^j\varphi_\ast\mathscr F)^{(-1)^j}=\det\varphi_\ast(\mathscr F)\otimes (\det R^1\varphi_\ast(\mathscr F))^{-1}
\end{equation}
Unfortunately equation (\ref{wrongdet_coho}) doesn't make any sense in general, since the higher direct images $R^j\varphi_\ast\mathscr F$ are not locally free sheaves, so we cannot take the determinant. However, we will cook up a definition of $\det R\varphi_\ast(\mathscr F)$ which agrees  with equation (\ref{wrongdet_coho}) when $R^j\varphi_\ast\mathscr F$ are locally free. 

The following  proposition is fundamental:
\begin{proposition}\label{fund_det_coho}
There exists  an effective divisor $D$ on $X$ which doesn't contain any fibre of $\varphi$ such that for any coherent sheaf $\mathscr F $ on $X$ we get an exact sequence:
\begin{equation}\label{fund_ex_seq}
0\to\varphi_\ast \mathscr F\to\varphi_\ast\mathscr F(D)\to \varphi_\ast(\mathscr F(D)/\mathscr F)\to R^1\varphi_\ast\mathscr F\to 0
\end{equation}
such that  $\varphi_\ast\mathscr F(D)$ and $\varphi_\ast(\mathscr F(D)/\mathscr F)$ are both locally free sheaves on $B$.
\end{proposition}
\proof
See  \cite[XIII section 4.]{algcu} or \cite[VI, Lemma 1.1]{Lang}.
\endproof
\begin{definition}
Let $\mathscr F\in\Coh(X)$ and let $D$ be a divisor as in proposition \ref{fund_det_coho}; \emph{the determinant of the cohomology of $\mathscr F$} is:
$$\det R\varphi_\ast(\mathscr F):=\det\varphi_\ast(\mathscr F(D))\otimes (\det \varphi_\ast(\mathscr F(D)/\mathscr F))^{-1}\in \Pic(B)$$
Moreover $\det R\varphi_\ast(\mathscr F)$ doesn't depend on the choice of $D$ (for the proof of this statement see \cite[XIII section 4.]{algcu} or \cite[VI]{Lang}).
\end{definition}

Now we want to show that if $R^0\varphi_\star \mathscr F=\varphi_\ast\mathscr F$ and  $R^1\varphi_\star \mathscr F$ are both locally free, then $\det R\varphi_\ast(\mathscr F)$ is given by equation (\ref{wrongdet_coho}). Consider the exact sequence of equation (\ref{fund_ex_seq}), put 
$$f:\varphi_\ast \mathscr F(D)\to\varphi_\ast(\mathscr F(D)/\mathscr F)\,,$$ 
$$g:\varphi_\ast(\mathscr F(D)/\mathscr F)\to R^1\varphi_\ast\mathscr F\,,$$
and $\mathscr G=\im(f)=\ker(g)$. Then we get the following two short exact sequences of locally free sheaves:
\begin{equation}\label{sh_ex_seq1}
0\to\varphi_\ast\mathscr F\to\varphi_\ast \mathscr F(D)\xrightarrow{f} \mathscr G\to 0\,;
\end{equation}
\begin{equation}\label{sh_ex_seq2}
0\to\mathscr G\to\varphi_\ast (\mathscr F(D)/\mathscr F)\xrightarrow{g} R^1\varphi_\ast\mathscr F\to 0\,.
\end{equation}
At this point we use the properties of the determinant on short exact sequences and we obtain:
$$\det R\varphi_\ast(\mathscr F)=\det\varphi_\ast(\mathscr F(D))\otimes (\det \varphi_\ast(\mathscr F(D)/\mathscr F))^{-1}\cong$$
$$\cong\det\varphi_\ast \mathscr F\otimes \det \mathscr G\otimes (\det \mathscr G)^{-1}\otimes (\det R^1\varphi_\ast(\mathscr F))^{-1}\cong$$
$$\cong \det\varphi_\ast \mathscr F\otimes \det R^1\varphi_\ast(\mathscr F))^{-1}\,.$$
The relationship between the determinant of cohomology and Deligne pairing is given by the following theorem:
\begin{theorem}\label{det_coh_deli}
Let $\mathscr L,\mathscr M$  be two invertible sheaves on $X$, then
$$\left<\mathscr L,\mathscr M\right>\cong \det R\varphi_\ast(\mathscr O_X)\otimes (\det R\varphi_\ast(\mathscr L))^{-1}\otimes(\det R\varphi_\ast(\mathscr M))^{-1}\otimes\det R\varphi_\ast(\mathscr L\otimes \mathscr M)\,.$$
\end{theorem}
\proof
See \cite[XIII, Theorem 5.8]{algcu}.
\endproof
\end{appendices}

\bibliographystyle{plain}
\bibliography{adI.bib}

\end{document}